\documentclass[sn-mathphys-num]{sn-jnl}
\usepackage{graphicx}
\usepackage{multirow}
\usepackage{amsmath,amssymb,amsfonts}
\usepackage{amsthm}
\usepackage{mathrsfs}
\usepackage[title]{appendix}
\usepackage{xcolor}
\usepackage{textcomp}
\usepackage{manyfoot}
\usepackage{booktabs}
\usepackage{algorithm}
\usepackage{algorithmicx}
\usepackage{algpseudocode}
\usepackage{listings}

\newtheorem{example}{Example}
\theoremstyle{case}

\usepackage{algorithm}
\usepackage{algpseudocode}
\usepackage{graphicx}
\usepackage{subcaption}
\usepackage{caption}

\theoremstyle{thmstyleone}

\theoremstyle{thmstyletwo}

\theoremstyle{thmstylethree}

\newtheorem{lemma}{Lemma}
\begin{document}

\title[{\bf A Variational Physics-Informed Neural Network Framework Using  Petrov-Galerkin Method for Solving Singularly Perturbed Boundary Value Problems}]{\bf A Variational Physics-Informed Neural Network Framework Using  Petrov-Galerkin Method for Solving Singularly Perturbed Boundary Value Problems}

\author[1]{\fnm{Vijay} \sur{Kumar}}\email{vijaykumarmathematics@gmail.com}

\author*[1]{\fnm{Gautam} \sur{Singh}}\email{gautam@nitt.edu}

\affil[1]{\orgdiv{Department of Mathematics}, \orgname{National Institute of Technology Tiruchirappalli}, \city{Tiruchirappalli}, \postcode{620015}, \state{Tamil Nadu}, \country{India}}

\abstract{ This work proposes a Variational Physics-Informed Neural Network (VPINN) framework that integrates the Petrov-Galerkin formulation with deep neural networks (DNNs) for solving one-dimensional singularly perturbed boundary value problems (BVPs) and parabolic partial differential equations (PDEs) involving one or two small parameters. The method adopts a nonlinear approximation in which the trial space is defined by neural network functions, while the test space is constructed from hat functions. The weak formulation is constructed using localized test functions, with interface penalty terms introduced to enhance numerical stability and accurately capture boundary layers. Dirichlet boundary conditions are imposed via hard constraints, and source terms are computed using automatic differentiation. Numerical experiments on benchmark problems demonstrate the effectiveness of the proposed method, showing significantly improved accuracy in both the $L_2$ and maximum norms compared to the standard VPINN approach for one-dimensional singularly perturbed differential equations (SPDEs).}

\keywords{ Variational physics-informed neural network, Singularly perturbed problem, Petrov-Galerkin method, Automatic differentiation.}

\pacs[MSC Classification]{34D15, 35B25, 65M12, 68T07.}

\maketitle

\section{Introduction}\label{secint}
Singular perturbation problems (SPPs) arise in a wide range of scientific and engineering disciplines, including chemical and nuclear engineering \cite{Alhu}, control theory \cite{kokotovic}, fluid flow modeling \cite{Modeling_Fluid}, water pollution studies \cite{Water_pollution}, and fluid dynamics \cite{Schlicting}. In SPPs, the presence of small parameters multiplying the highest-order derivatives leads to boundary layer behavior, characterized by sharp variations near the boundaries or within the interior of the domain.
This makes solving such problems particularly challenging.

Many numerical techniques have been developed over the years to solve one- and two-parameter second-order singularly perturbed BVPs, particularly using the finite element method (FEM) \cite{Parabolic_FEM_2D,Melenk,Roos2, Gautam1, Gautam4, fem_tobiska,FEM_2D}. Among these, the Discontinuous Galerkin (DG) method introduced by Reed and Hill in 1973 to address neutron transport problems \cite{Reed} has emerged as a highly flexible and powerful technique for solving various classes of PDEs. To address the challenges posed by singular perturbations, several DG variants have been proposed, such as the Interior Penalty Discontinuous Galerkin (IPDG) method \cite{IPDG}, the Symmetric Interior Penalty Galerkin (SIPG) method \cite{SIPG}, and the Non-Symmetric Interior Penalty Galerkin (NIPG) method \cite{Gautam_parabolic}. For further details on the  DG method, readers may refer to book \cite{Riviere}. These methods are particularly effective in handling discontinuities and capturing steep solution gradients, making them well-suited for problems exhibiting boundary or interior layers.

In parallel, Deep Neural Networks (DNNs) have emerged as powerful function approximators capable of learning continuous solutions directly from governing equations. This paradigm shift led to the development of Physics-Informed Neural Networks (PINNs), introduced by Raissi et al. \cite{PINN_Raissi}, where neural networks are trained to satisfy the differential equations and associated boundary or initial conditions by minimizing a physics-informed loss function. As a mesh-free approach, PINNs can naturally accommodate irregular domains and high-dimensional problems, making them appealing for many applications in scientific computing.

The Finite Basis Physics-Informed Neural Network (FB-PINN) framework \cite{FB_PINN} incorporates finite element methodology by approximating the solution as a sum of finite basis functions. This structure enhances PINN performance, particularly for multiscale problems, and was extended by Raina et al. \cite{FB-PINN_Natesan} to handle complex two-parameter convection-dominated systems. In parallel, to address the limitations of standard PINNs in capturing sharp boundary layers in SPPs, Cao et al. \cite{PAPINN} introduced the Parameter Asymptotic PINN (PA-PINN) for one-parameter cases. This approach was later generalized by Boro et al. \cite{PINN_Natesan} to effectively solve two-parameter convection-dominated problems, demonstrating improved accuracy through the incorporation of asymptotic information.

Moreover, for problems involving non-smooth solutions, several PINN variants based on the variational (weak) formulation have been proposed. These include a Petrov–Galerkin method that employs DNNs as nonlinear trial functions and uses Legendre polynomials or trigonometric functions as test functions \cite{VPINN1}. The hp-VPINN framework \cite{VPINN} introduces  hp-refinement by decomposing the computational domain and employing projections onto high-order Legendre polynomial spaces. It constructs the trial space using a globally defined neural network, while the test space consists of locally supported piecewise polynomial functions on each subdomain.

Yadav and Ganesan \cite{Yadav2} proposed an artificial neural network (ANNs) augmented stabilization strategy for the Streamline Upwind/Petrov–Galerkin FEM to effectively solve SPPs. Chen et al. \cite{PINN_NIPG} introduced a Discontinuous Galerkin-induced Neural Network (DGNN). However, to the best of our knowledge, no article has been published on solving SPPs using the VPINN framework with the Petrov-Galerkin method

In this paper, we develop a VPINN with the Petrov-Galerkin framework to solve one-parameter and two-parameter singularly perturbed BVPs, as well as parabolic PDEs involving small perturbation parameters. The method employs a nonlinear approximation by selecting the trial space as the space of neural networks and the test space as the space of hat functions. Dirichlet boundary conditions are imposed via hard constraints, and source terms are computed through automatic differentiation. For time-dependent parabolic problems, temporal discretization is performed using the backward Euler scheme.

The structure of the paper is as follows: Section \ref{sec1} provides an overview of ANNs, PINNs, VPINNs, and the dataset associated with the problem under study. Section \ref{sec3} introduces the model problems, the proposed algorithm, and numerical experiments. Finally, Section \ref{Conclusion} concludes the paper.

\section{Basics of neural network frameworks}\label{sec1}

\subsection{Artificial neural networks (ANNs)}\label{ANN}
Let $F:\mathbb{R}^m \to \mathbb{R}^{mout}$ denotes a feed-forward neural network with $s-1$ hidden layers and $nL$ neurons in the  $L$-th layer. Usually, the first and last layers are referred to as the input and output layers, respectively, while all intermediate layers are known as hidden layers. When a network has more than one hidden layer, it is commonly referred to as a DNN. ANNs are widely recognized for their capacity to approximate multivariate functions and are often employed as universal function approximators. In a typical fully connected feedforward network, an input \( x \in \mathbb{R}^d \) is mapped to an output \( u_\theta(x) \in \mathbb{R} \) via a sequence of affine transformations and nonlinear activation functions. For a network with $L$ layers, the forward pass is defined as:
\[
u_\theta(x) = \mathcal{N}_L \circ \Phi \circ \mathcal{N}_{L-1} \circ \cdots \circ \Phi \circ \mathcal{N}_1(x),
\]
\[
\mathcal{N}_k(x) = W_k \Phi(\mathcal{N}_{k-1}(x)) + b_k, \quad \mbox{for}\quad  2 \leq L \leq s \quad \text{with } \mathcal{N}_0(x) = x,
\]
where $W_k \in \mathbb{R}^{N_k \times N_{k-1}}$ and $b_k \in \mathbb{R}^{N_k}$ denote the trainable weight matrix and bias vector corresponding to layer  $k$. The activation function \( \Phi \) is a nonlinear map, such as ReLU, sigmoid, hyperbolic tangent (tanh), or leaky ReLU. The collection of all trainable parameters is denoted by \( \theta = \{W_k, b_k\}_{k=1}^{L} \).

In traditional supervised learning, the neural network parameters are optimized by minimizing a loss function over a dataset of labeled examples. However, in the context of scientific machine learning where labeled data may be limited an alternative approach is to encode known physical laws directly into the loss function. A Multilayer Perceptron (MLP) is a feedforward neural network that connects every neuron in one layer to all neurons in the next layer.

\begin{figure}[htbp]
	\centering
	\includegraphics[width=0.7\textwidth]{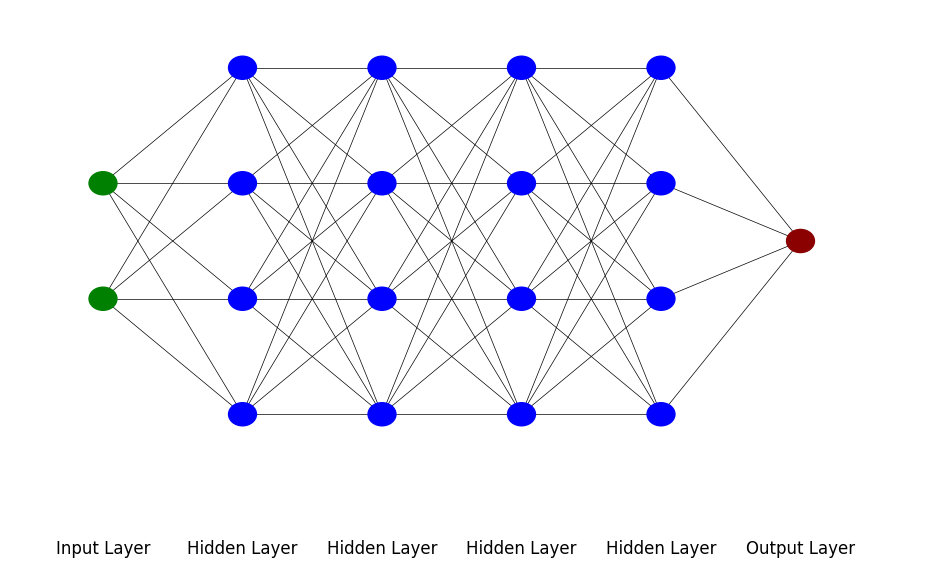}
	\label{fig:ANNs}
	\caption{Visualization of artificial neural network.}\label{figureANNs}
\end{figure}

ANNs have been widely recognized as powerful tools for addressing complex and unresolved challenges in scientific computing, owing to their universal approximation capabilities and the flexibility offered by their hyper-parameter tuning. ANNs have been effectively employed for solving PDEs \cite{DNN_FEM_1, DNN_FEM_2}. Recent studies have also investigated hybrid approaches that integrate ANNs with classical numerical techniques \cite{DNN_FEM_3}. For instance, one study proposed the Proper orthogonal decomposition neural network  (POD-NN) method, where ANNs are utilized during the interpolation stage to develop a non-intrusive Radial Basis (RB) method based on proper orthogonal decomposition. Similarly, \cite{DNN_FEM_4} introduced a smart FEM that learns a direct mapping between element-level inputs and outputs to alleviate the computational burden of large-scale problems. However, in both approaches, FEM is not directly embedded within the deep learning framework, but rather employed for data generation or as a means of accelerating specific numerical operations. A visualization of an ANNs is provided in Figure~\ref{figureANNs}.

\subsection{Physics-informed neural networks (PINNs)}\label{PINN}

PINNs, as introduced in \cite{PINN_Raissi}, provide a robust approach for addressing both forward and inverse problems involving differential and integro-differential equations. They are particularly effective in scenarios where data is limited, uncertain, or derived from multiple sources of varying fidelity. One of the defining advantages of PINNs is their ability to integrate diverse sources of information such as governing equations, boundary and initial conditions, and empirical observations directly into the model's loss function. This integration constrains the solution space, ensuring that the learned solution adheres to the underlying physical laws. Since PINNs are mesh-free, they recast the task of solving PDEs as an optimization problem over a loss function, rather than relying on traditional discretization techniques. As a result, they can achieve strong generalization with fewer data points compared to conventional purely data-driven neural networks.

Despite their effectiveness, applying PINNs to SPPs remains challenging. This is due to the nature of SPPs, where the solution typically exhibits rapid variations near boundary or interior layers, while varying more smoothly elsewhere. Capturing such multi-scale behavior with a single deep neural network can be difficult. A practical approach in such cases is to divide the computational domain into multiple regions, with each region modeled independently based on its specific characteristics. This strategy can improve the solution’s accuracy and efficiency. Nevertheless, domain decomposition techniques, such as those proposed in \cite{Deep_decomposition}, often depend on prior knowledge of the layer structure or singularity locations information that is not always available in practical applications.

PINNs incorporate prior knowledge of physical laws such as partial differential equations (PDEs) into the neural network training. Consider a PDE of the form:
\[
\mathcal{D}(u(x)) = r(x), \quad x \in I \subset \mathbb{R}^d, \qquad \mathcal{B}(u(x)) = 0, \quad x \in \partial I,
\]
where \( \mathcal{D} \) is a differential operator and \( \mathcal{B} \) denotes the boundary conditions. The neural network \( u_\theta(x) \) is used to approximate the solution. The residual is defined by:
\[
v(x) := \mathcal{D}(u_\theta(x)) - r(x).
\]

The total PINN loss function includes two components:
\[
\text{Loss}_{\text{PINN}}(\theta) = \frac{1}{N_f} \sum_{i=1}^{N_f} \left| v(x_i^f) \right|^2 + \frac{1}{N_u} \sum_{i=1}^{N_u} \left| u_\theta(x_i^u) - g(x_i^u) \right|^2,
\]
where \( \{x_i^f\}_{i=1}^{N_f} \subset I \) are the \emph{interior collocation points} at which the differential equation is enforced, and \( \{x_i^u\}_{i=1}^{N_u} \subset \partial I \) are the \emph{boundary points} where the Dirichlet boundary condition \( u(x) = g(x) \) is imposed. Here, \( N_f \) and \( N_u \) denote the total number of interior and boundary points, respectively.

To improve accuracy at the boundary, one may enforce hard constraints by modifying the neural network architecture so that the boundary conditions are exactly satisfied, i.e., by constructing \( u_\theta(x) \) such that \( u_\theta|_{\partial I} = g(x) \). This eliminates the need to include a penalty term for boundary conditions in the loss and can improve convergence and accuracy.

\subsection{Variational physics-informed neural networks (VPINNs)}\label{VPINN}

VPINNs extend the PINN framework by enforcing the residual in a variational or weak form. Instead of minimizing $u(x)$ pointwise, VPINNs require that the residual be orthogonal to a set of test functions $\{v_j(x)\}_{j=1}^M$, leading to:
\[
\int_{I} u(x) \, v_j(x) \, dx = 0 \quad \forall  j = 1, \ldots, M,
\]
that is,
\[
\int_{I} \left( \mathcal{D}(u_\theta(x)) - r(x) \right) v_j(x) \, dx = 0.
\]

The corresponding VPINN loss function is:
\[
\text{Loss}_{\text{VPINN}}(\theta) = \sum_{j=1}^{M} \left| \int_{I} \left( \mathcal{D}(u_\theta(x)) - r(x) \right) v_j(x) \, dx \right|^2 + \text{MSE}_u.
\]

As in PINNs, the boundary loss term is:
\[
\text{MSE}_u = \frac{1}{N_u} \sum_{i=1}^{N_u} \left| u_\theta(x_i^u) - g(x_i^u) \right|^2.
\]
Alternatively, using hard constraint boundary conditions can be implemented in VPINNs as well by ensuring the network structure satisfies $u_\theta(x) = g(x)$ on $\partial I$ exactly.

The test functions $v_j$ are typically selected from a piecewise polynomial basis, such as linear hat functions or higher-order finite element bases. All integrals are evaluated using suitable numerical quadrature schemes.

\subsection{Dataset}\label{data}
This subsection focuses on data generated from one-dimensional singularly perturbed differential equations. In contrast to other methods, our VPINN framework employs $36$ test functions and $1000$ quadrature points. A hard boundary constraint is applied to enforce Dirichlet conditions. For every neural network architecture featuring $4$ hidden layers with $20$ neurons each and the $\tanh$ activation function. For the convection-diffusion problem, the model is trained using the Adam optimizer with a learning rate of  $10^{-3}$ and is then fine-tuned using the L-BFGS optimizer for $1500$ epochs. In contrast, for the reaction-diffusion and two-parameter singularly perturbed problems, only the L-BFGS optimizer is used for $1500$ epochs.

For time-dependent problems, the framework employs $18$ test functions and incorporates $100$ quadrature points in both space and time. Training begins with the Adam optimizer employing a learning rate of  $10^{-3}$ and is then fine-tuned using the L-BFGS optimizer for $1000$ epochs, while maintaining the same network architecture and dataset structure as in the static case.

\section{\textbf{ Model problems and numerical experiments}}\label{sec3}

\subsection{The convection-diffusion BVPs}\label{sec3.1}
We examine the following one-parameter singularly perturbed BVPs:
\begin{equation}\label{E1}
	\left\{
	\begin{array} {ll}
		-\epsilon u''(x)+b(x)u'(x)+c(x)u(x)=r(x),\  x\in I=(0,\,1), \\[8pt]
		u(0)=u(1)=0.\\[8pt]
	\end{array}\right.
\end{equation}
Here, $0<\epsilon\ll 1$ is a small parameter, and the function  $b(x), c(x)$ and  $r(x)$ are sufficiently smooth and satisfy the conditions  $b(x)\geq \alpha > 0$, \quad $c(x)-\frac{b'(x)}{2}>0$,  $\forall x \in I$, for constant $\alpha$. Therefore, the solution $u(x)$ of (\ref{E1}) features a boundary layer at  $x=1$.

The finite element formulation of the BVPs~\eqref{E1} seeks a function \( u \in H_0^1(I) \) such that, $\forall$ \( v \in H_0^1(I) \),
\begin{equation}\label{E2}
	a(u, v) = (r, v),
\end{equation}
where the bilinear form \( a(\cdot, \cdot) \) is given by
\[
a(u, v) = \int_{I} \epsilon\, u'(x)\, v'(x)\, dx 
+ \int_{I} b(x)\, u'(x)\, v(x)\, dx 
+ \int_{I} c(x)\, u(x)\, v(x)\, dx,
\]
and the linear functional on the right-hand side is defined as
\[
(r, v) = \int_{I} r(x)\, v(x)\, dx.
\]

We denote by \( I_h \) a mesh that subdivides the domain \( I \) into a finite number of elements. Corresponding to this mesh, we define a finite-dimensional subspace \( V_h \subset H_0^1(I) \). We also define the neural network-based trial space as $U_{NN} = \left\{ u_\theta(x) = \mathrm{DNN}(x; \theta) \,\middle|\, \theta \in \Theta \right\},$
where \( \mathrm{DNN}(x; \theta) \) denotes a deep neural network formed by sequential layers of linear mappings combined with nonlinear activation functions,  and \(  \Theta \) represents the collection of all trainable parameters.

The discrete variational problem associated with equation~\eqref{E2} is then stated as follows:

Find \( u_h \in U_{NN} \) such that, $\forall$ \( v \in V_h \),
\[
a_h(u_h, v) = (r, v),
\]
where the notation \( a_h(\cdot, \cdot) \) implies that integration is performed element-wise over the partition \( I_h \).

\begin{figure}[htbp]
	\centering
	\includegraphics[width=1\textwidth]{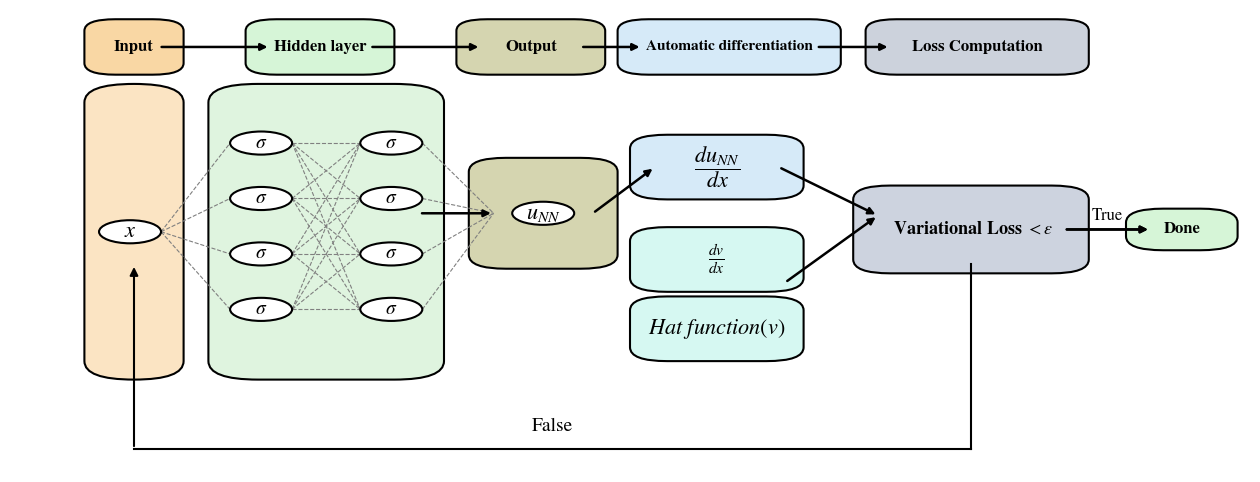}
	\caption{VPINN algorithm flowchart.}\label{figure_algorithm1}
\end{figure}

The VPINN algorithm for one-parameter  convection-diffusion problem is illustrated in Algorithm \ref{Algorithm1}. Figure \ref{figure_algorithm1} provides a visual representation of the  Algorithm \ref{Algorithm1}.

\begin{algorithm}
	\caption{VPINN Algorithm for Solving $-\epsilon u''(x) +b(x)u'(x)+ c(x)u(x) = r(x)$ in $(0,1)$ with $u(0) = u(1) = 0$}
	\begin{algorithmic}[1]
		\State \textbf{Input:} Small parameter $\epsilon > 0$, number of collocation points $N$, number of test functions $M$
		\State \textbf{Define:} Trial function $\mathcal{T}_\theta(x) := x(1 - x)\cdot \text{NN}_\theta(x)$ to satisfy boundary conditions
		\State \textbf{Generate:} Uniformly spaced collocation points $\{x_j\}_{j=1}^{N} \subset (0,1)$
		\State \textbf{Construct:} Piecewise linear hat test functions $\{v_i(x)\}_{i=1}^{M}$ on a mesh $\{x_k\}_{k=0}^{M+1}$
		\State \textbf{Compute:} Exact solution $u_{\text{exact}}(x)$ and source term $r(x) = -\epsilon u_{\text{exact}}''(x) +b(x) u_{\text{exact}}'(x) +c(x)u_{\text{exact}}(x)$ using automatic differentiation
		\For{$i = 1$ to $M$}
		\State Define weak residual for each test function $v_i(x)$:
		\[
		R_i(\theta) := \int_0^1 \left(-\epsilon^2 \mathcal{T}_\theta''(x) +b(x) \mathcal{T}_\theta'(x)+ c(x)\mathcal{T}_\theta(x) - r(x)\right)v_i(x)\,dx
		\]
		\EndFor
		\State \textbf{Define:} VPINN loss function as:
		\[
		\mathcal{L}(\theta) := \frac{1}{M} \sum_{i=1}^M \left( R_i(\theta) \right)^2
		\]
		\State \textbf{Initialize:} Neural network weights $\theta$ randomly
		\State \textbf{Train:} Minimize $\mathcal{L}(\theta)$ first using Adam and then fine-tuned using the L-BFGS optimizer with closure-based backpropagation
		\State \textbf{Return:} Trained network $\mathcal{T}_\theta(x)$ approximating the solution $u(x)$
	\end{algorithmic}\label{Algorithm1}
\end{algorithm}

\begin{example}\label{example1}
Consider the problem in \eqref{E1} with
\( b(x) = -(1+x) \) and \( c(x) = 1 \),
where the source function \( r(x) \) is computed from the exact solution
\[
u(x) = 1 - x - e^{-x/\epsilon} + x\,e^{-x/\epsilon}.
\]
\end{example}

\begin{figure}[htbp]
\centering

\begin{subfigure}[b]{0.3\textwidth}
	\includegraphics[width=\textwidth,height=1.2\textwidth]{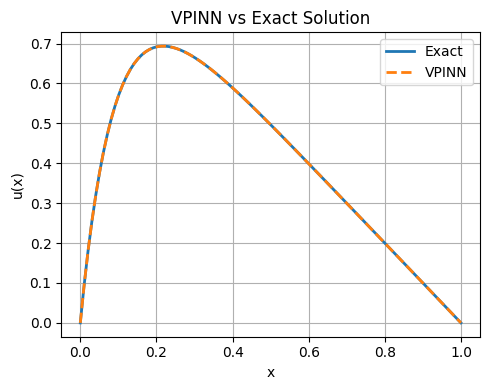}
	\caption{}
	\label{fig:subfig1A}
\end{subfigure}
\begin{subfigure}[b]{0.3\textwidth}
	\includegraphics[width=\textwidth,height=1.2\textwidth]{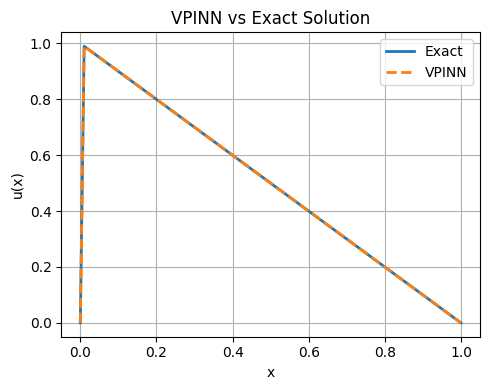}
	\caption{}
	\label{fig:subfig1B}
\end{subfigure}
\begin{subfigure}[b]{0.3\textwidth}
	\includegraphics[width=\textwidth,height=1.2\textwidth]{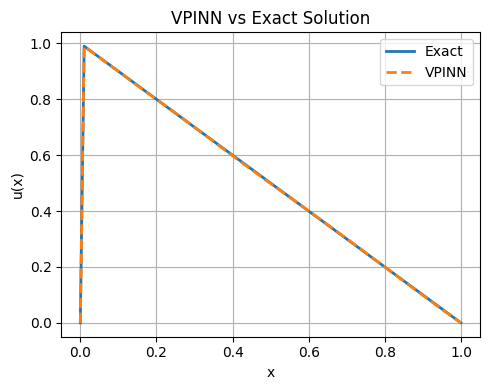}
	\caption{}
	\label{fig:subfig1C}
\end{subfigure}
\caption{Comparison between Exact and VPINN solution for Example \ref{example1} with different value of perturbation parameter $\epsilon=10^{-1},10^{-3}$, and $10^{-5}$, corresponding to subfigures~(a), (b), and (c), respectively.}\label{figure1}
\end{figure}

\begin{figure}[htbp]
\centering

\begin{subfigure}[b]{0.3\textwidth}
	\includegraphics[width=\textwidth,height=1.0\textwidth]{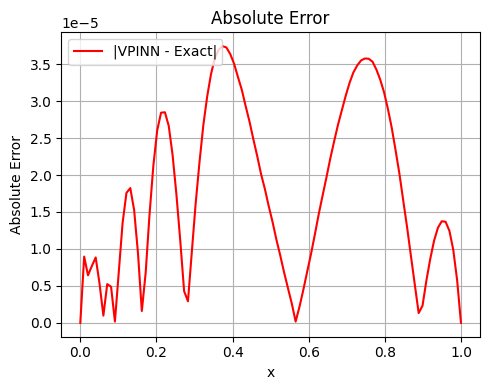}
	\caption{}
	\label{fig:subfig2A}
\end{subfigure}
\hfill
\begin{subfigure}[b]{0.3\textwidth}
	\includegraphics[width=\textwidth,height=1.0\textwidth]{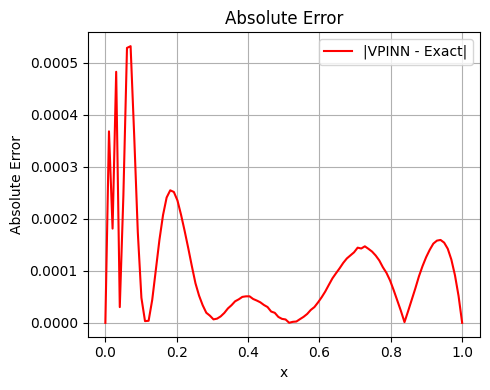}
	\caption{}
	\label{fig:subfig2B}
\end{subfigure}
\hfill
\begin{subfigure}[b]{0.3\textwidth}
	\includegraphics[width=\textwidth,height=1.0\textwidth]{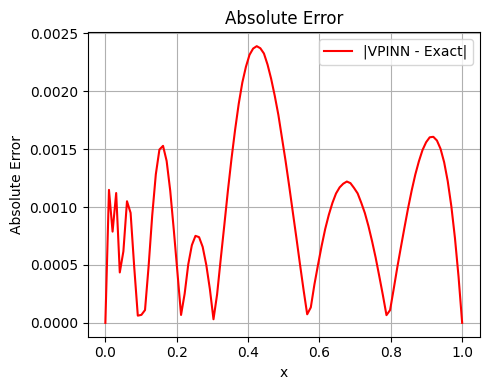}
	\caption{}
	\label{fig:subfig2C}
\end{subfigure}
\caption{Error vs x for Example \ref{example1} with different value of perturbation parameter $\epsilon=10^{-1},10^{-3}$, and $10^{-5}$, corresponding to subfigures~(a), (b), and (c), respectively.}\label{figure2}
\end{figure}

\begin{figure}[htbp]
\centering

\begin{subfigure}[b]{0.3\textwidth}
	\includegraphics[width=\textwidth,height=1.0\textwidth]{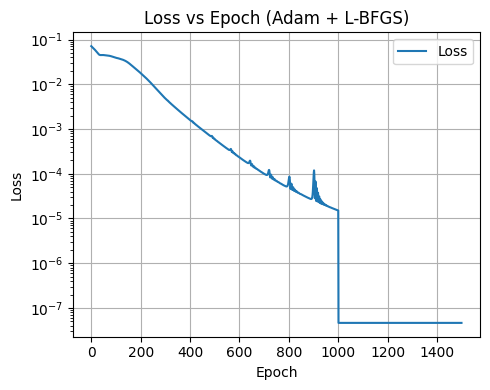}
	\caption{}
	\label{fig:subfig3A}
\end{subfigure}
\hfill
\begin{subfigure}[b]{0.3\textwidth}
	\includegraphics[width=\textwidth,height=1.0\textwidth]{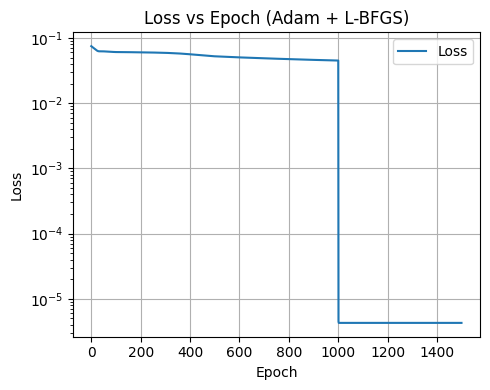}
	\caption{}
	\label{fig:subfig3B}
\end{subfigure}
\hfill
\begin{subfigure}[b]{0.3\textwidth}
	\includegraphics[width=\textwidth,height=1.0\textwidth]{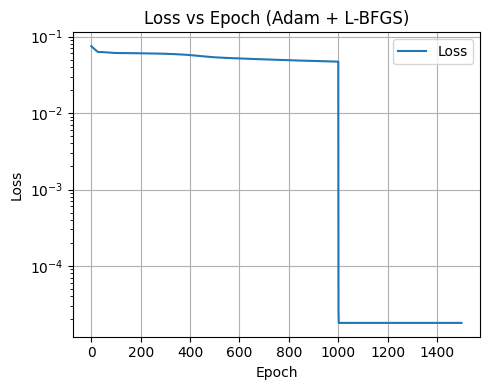}
	\caption{}
	\label{fig:subfig3C}
\end{subfigure}
\caption{Loss vs Epoch for Example \ref{example1} with different value of perturbation parameter $\epsilon=10^{-1},10^{-3}$, $10^{-5}$, corresponding to subfigures~(a), (b), and (c), respectively.}\label{figure3}
\end{figure}

We executed the code for $1500$ iterations and the Maximum error, Relative maximum error, $L_{2}$ error, Relative $L_{2}$ error and Loss at $1500$ epoch obtained are presented in Table \ref{tbl_1}. Figure \ref{figure1} presents a comparison between the exact solution and VPINN approximation for different parameter values. Figure \ref{figure2} illustrates the absolute error plots for the same example, while Figure \ref{figure3} shows the loss versus epoch curves, corresponding to the final iteration of our algorithm.

\begin{table}[htbp]
\caption{\it{ Loss and Error for Example \ref{example1}.}}\label{tbl_1}  
\begin{tabular}{@{}lllllll@{}}
	\multicolumn 1 {c}{}  & \multicolumn 5 {c}
	{}\\
	\hline
	$\epsilon$	&Loss & maximum & Rel-maximum  & $L_{2}$ & Rel-$L_{2}$ \\
	\hline
	$10^{-1}$ & 4.64323e-08  & 3.74317e-05  & 5.39757e-05 & 2.11394e-05 & 4.58645e-06 \\ [6pt]
	
	$10^{-2}$  &  3.82179e-07 &2.89887e-04  & 3.07275e-04 & 1.43689e-04 & 2.55847e-05 \\ [6pt]
	
	$10^{-3}$   & 4.27739e-06 & 5.31852e-04   & 5.37301e-04 & 1.45876e-04 & 2.55877e-05   \\ [6pt]
	
	$10^{-4}$  & 1.64225e-05 & 1.77413e-03  &  1.79224e-03 & 9.12843e-04  &  1.60119e-04 \\ [6pt]

	$10^{-5}$  &1.80347e-05   &2.38985e-03  & 2.41423e-03 & 1.16507e-03 & 2.04361e-04  \\ [6pt]
	
	\hline
\end{tabular}
\end{table}

\subsection{Initial boundary-value problem (IBVP)} \label{sec3.2}
Consider the following singularly perturbed parabolic PDEs with
one-parameters in the domain $G$ = $I$ $\times$ $(0, T]$, where $I = (0, 1):$
\begin{equation}\label{E4}
\left\{
\begin{array} {ll}
	u_{t}	-\epsilon u_{xx}+b(x)u_{x}+c(x)u(x,t)=r(x,t),\  x\in I=(0,\,1), \\[8pt]
	u(0,t)= u(1,t)=0,  \quad 0\leq t \leq T,\\[8pt]
	u(x,0)-\phi(x)=0, \quad \mbox{for}\quad  x \in [0, 1],
\end{array}\right.
\end{equation}
Here, $0<\epsilon\ll 1$ is a small parameter, and the function  $b(x), c(x)$ and  $r(x,t)$ are sufficiently smooth and satisfy the conditions  $b(x)\geq \alpha > 0$, \quad $c(x)-\frac{b'(x)}{2}>0$,  $\forall x \in I$, for constant $\alpha$. Therefore solution $u(x,t)$ of (\ref{E4}) exhibits boundary layers of at $x=1$.

The finite element formulation of the BVPs (\ref{E4}) is given by. Find  $u\in H_0^1(I) $, such that $\forall v \in  H_0^1(I)$ 
\begin{equation}\label{E5}
a(u,v)=(r,v),
\end{equation}
where the bilinear form $a(.,.)$ is defined by
\[ \int_{I}\frac{\partial{u}}{\partial{t}}vdx+\int_{I}\epsilon\frac{\partial{u}}{\partial{x}}\frac{\partial{v}}{\partial{x}}dx+\int_{I}b(x)\frac{\partial{u}}{\partial{x}}vdx+\int_{I}c(x)uvdx=\int_{I}r(x)vdx,\]
At each time step $t=t_{m}$ $(1\leq m\leq M)$, we have the following
\begin{equation}\label{E6}
a(u^{m},v)=(r^{m},v),
\end{equation}
\[ \int_{I}\frac{\partial{u^{m}}}{\partial{t}}vdx+\int_{I}\epsilon\frac{\partial{u^{m}}}{\partial{x}}\frac{\partial{v}}{\partial{x}}dx+\int_{I}b(x)\frac{\partial{u^{m}}}{\partial{x}}vdx+\int_{I}c(x)u^{m}vdx=\int_{I}r^{m}vdx,\]
\[ \int_{I}\frac{(u^{m}-u^{m-1})}{\Delta t}vdx+\int_{I}\epsilon\frac{\partial{u^{m}}}{\partial{x}}\frac{\partial{v}}{\partial{x}}dx+\int_{I}b(x)\frac{\partial{u^{m}}}{\partial{x}}vdx+\int_{I}c(x)u^{m}vdx=\int_{I}r^{m}vdx,\]
\[b(u^{m},v)=(g^{m},v)\]
\[b(u^{m},v)=\int_{I}\epsilon\frac{\partial{u^{m}}}{\partial{x}}\frac{\partial{v}}{\partial{x}}dx+\int_{I}b(x)\frac{\partial{u^{m}}}{\partial{x}}vdx+\int_{I}\bigg(c(x)+\frac{1}{\Delta t}\bigg)u^{m}vdx\]
\[(g^{m},v)=\int_{I}\bigg(r^{m}+\frac{u^{m-1}}{\Delta t}\bigg)v(x)dx.\]

We denote by \( I_h \) a mesh that subdivides the domain \( I \) into a finite number of elements. The discrete variational problem associated with equation~\eqref{E6} is stated as follows:

Find \( u_h \in U_{NN} \) such that, $\forall$ \( v \in V_h \),
\[ b_{h}(u^{m}_{h},v)=(g^{m},v).\]
where the notation \( b_h(\cdot, \cdot) \) implies that integration is performed element-wise over the partition \( I_h \).

The VPINN algorithm for Parabolic convection-diffusion IBVP is illustrated in Algorithm \ref{Algorithm2}.

\begin{algorithm}
\caption{VPINN for Solving $u_t - \varepsilon u_{xx} + b(x)u_x + c(x)u = r(x,t)$ with Dirichlet Boundary Conditions}
\label{alg:vpinn-parabolic}
\begin{algorithmic}[1]
	\State \textbf{Given:} PDE $u_t - \varepsilon u_{xx} + b(x)u_x + c(x)u = r(x,t)$ on $(x,t) \in (0,1)\times(0,1]$
	\State \textbf{Inputs:} Diffusion parameter $\varepsilon > 0$, spatial nodes $N_x$, time steps $N_t$, test functions $M$
	\State Discretize: $x_i \in [0,1]$, $t_n = n \Delta t$, with $\Delta t = \frac{1}{N_t}$
	\State Initialize solution $u^0(x) = u_{\text{exact}}(x,0)$
	\For{$n = 1$ to $N_t$}
	\State Define neural network trial solution $u_\theta(x) = x(1 - x)\cdot \text{NN}(x)$ to boundary condition.
	\State Compute source term $r(x, t_n)$ using autograd from $u_{\text{exact}}(x,t)$
	\State Construct residual:
	\[
	R(x) = \frac{u_\theta(x) - u^{n-1}(x)}{\Delta t} - \varepsilon u_{\theta,xx} + b(x)u_{\theta,x} +c(x) u_\theta(x) - r(x, t_n)
	\]
	\State Define weak loss using hat test functions $\{\phi_i(x)\}_{i=1}^M$:
	\[
	\mathcal{L}(\theta) = \frac{1}{M} \sum_{i=1}^M \int_0^1 R(x) \phi_i(x) \, dx \approx \frac{1}{M} \sum_{i=1}^M \text{MSE}(R(x)\phi_i(x))
	\]
	\Statex Mean Squared Error (MSE):
	\[
	\text{MSE}(g(x)) = \frac{1}{N_x} \sum_{j=1}^{N_x} \left(g(x_j)\right)^2
	\]
	\State Train the neural network by minimizing $\mathcal{L}(\theta)$ first using Adam and then fine-tuned using the L-BFGS optimizer
	\State Set $u^n(x) \gets u_\theta(x)$
	\EndFor
	\State \textbf{Return:} VPINN predictions $u_\theta(x,t_n)$ for $n = 1, \dots, N_t$
	\State \textbf{Postprocess:} Compute relative $L^2$ and $L^\infty$ errors, visualize solution surfaces
\end{algorithmic}\label{Algorithm2}
\end{algorithm}

%

\begin{example}\label{example2}
Consider  the parabolic one-parameter IBVP \eqref{E4} with 
\( b(x) = 1 \) and \( c(x) = 1 \).
The source term \( r(x,t) \) and initial condition \( \phi(x) \) 
are chosen from the exact solution
\[
u(x,t) = e^{-t}\big(1 - e^{-(1-x)/\epsilon}\big)\sin x .
\]
\end{example}

\begin{figure}[htbp]
\centering

\begin{subfigure}[b]{0.3\textwidth}
	\includegraphics[width=\textwidth,height=1.2\textwidth]{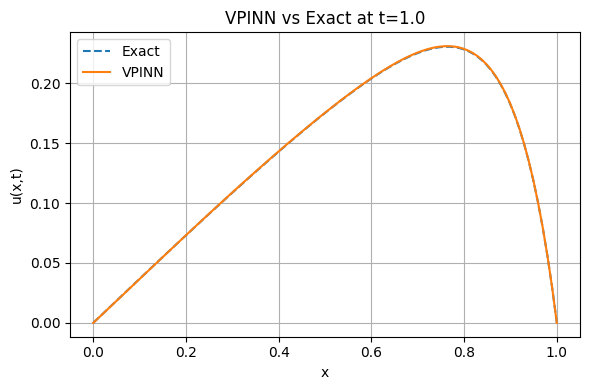}
	\caption{}
	\label{fig:subfig4A}
\end{subfigure}
\hfill
\begin{subfigure}[b]{0.3\textwidth}
	\includegraphics[width=\textwidth,height=1.2\textwidth]{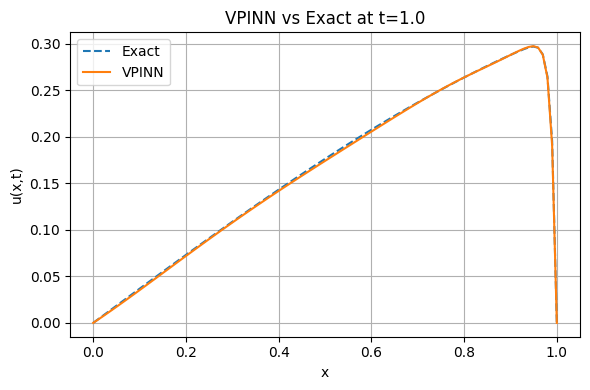}
	\caption{}
	\label{fig:subfig4B}
\end{subfigure}
\hfill
\begin{subfigure}[b]{0.3\textwidth}
	\includegraphics[width=\textwidth,height=1.2\textwidth]{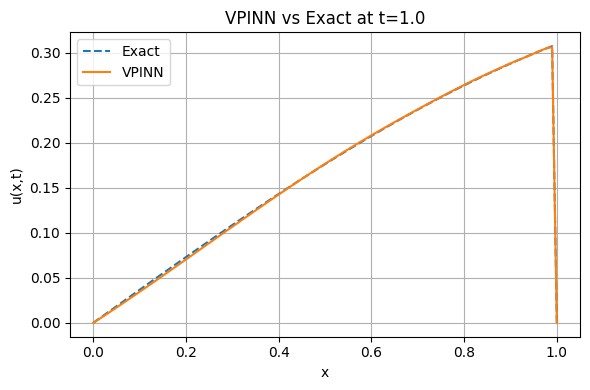}
	\caption{}
	\label{fig:subfig4C}
\end{subfigure}
\caption{Comparison between Exact and VPINN solution for Example \ref{example2} with different value of perturbation parameter $\epsilon=10^{-1},10^{-2}$, and $10^{-3}$, corresponding to subfigures~(a), (b), and (c), respectively.}\label{figure4}
\end{figure}

\begin{figure}[htbp]
\centering

\begin{subfigure}[b]{0.3\textwidth}
	\includegraphics[width=\textwidth,height=1.0\textwidth]{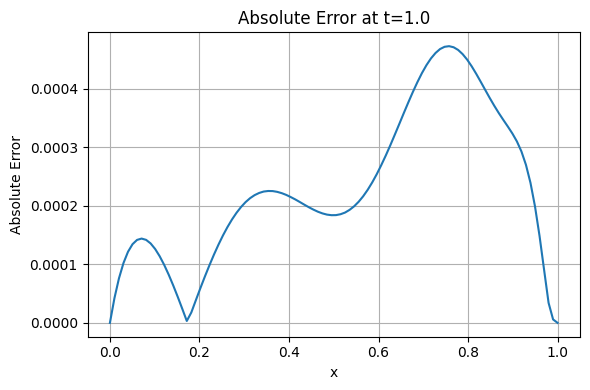}
	\caption{}
	\label{fig:subfig5A}
\end{subfigure}
\hfill
\begin{subfigure}[b]{0.3\textwidth}
	\includegraphics[width=\textwidth,height=1.0\textwidth]{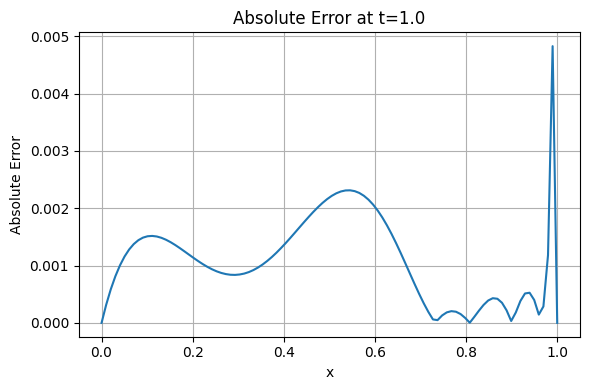}
	\caption{}
	\label{fig:subfig5B}
\end{subfigure}
\hfill
\begin{subfigure}[b]{0.3\textwidth}
	\includegraphics[width=\textwidth,height=1.0\textwidth]{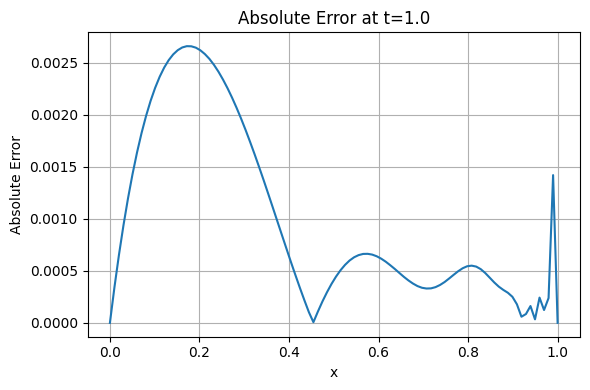}
	\caption{}
	\label{fig:subfig5C}
\end{subfigure}
\caption{Error vs x  for Example \ref{example2} with different value of perturbation parameter $\epsilon=10^{-1},10^{-2}$, and $10^{-3}$, corresponding to subfigures~(a), (b), and (c), respectively.}\label{figure5}
\end{figure}

\begin{figure}[htbp]
\centering

\begin{subfigure}[b]{0.3\textwidth}
	\includegraphics[width=\textwidth,height=1.0\textwidth]{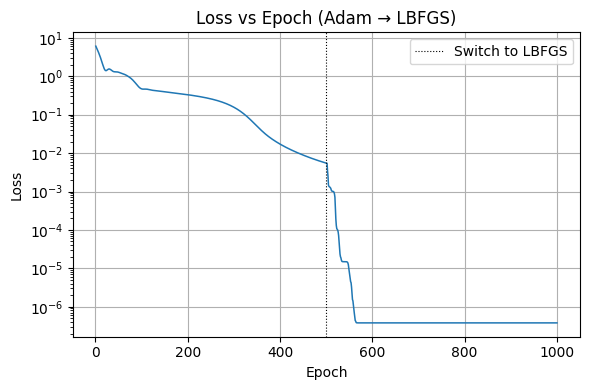}
	\caption{}
	\label{fig:subfig6A}
\end{subfigure}
\hfill
\begin{subfigure}[b]{0.3\textwidth}
	\includegraphics[width=\textwidth,height=1.0\textwidth]{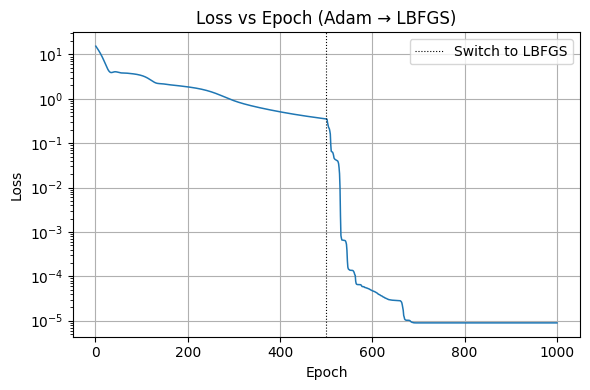}
	\caption{}
	\label{fig:subfig6B}
\end{subfigure}
\hfill
\begin{subfigure}[b]{0.3\textwidth}
	\includegraphics[width=\textwidth,height=1.0\textwidth]{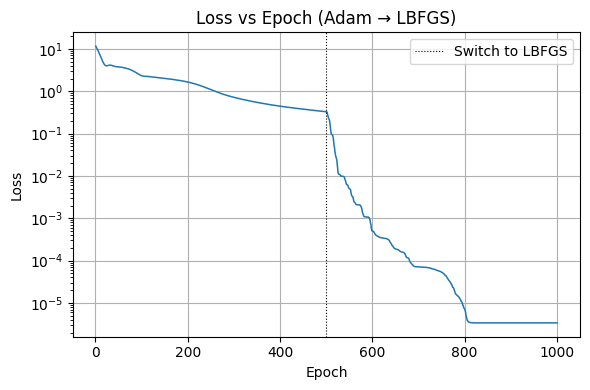}
	\caption{}
	\label{fig:subfig6C}
\end{subfigure}
\caption{Loss vs Epoch for Example \ref{example2} with different value of perturbation parameter $\epsilon=10^{-1},10^{-2}$,and $10^{-3}$, corresponding to subfigures~(a), (b), and (c), respectively.}\label{figure6}
\end{figure}

\begin{figure}[htbp]
\centering

\begin{subfigure}[b]{1.0\textwidth}
	\includegraphics[width=\textwidth]{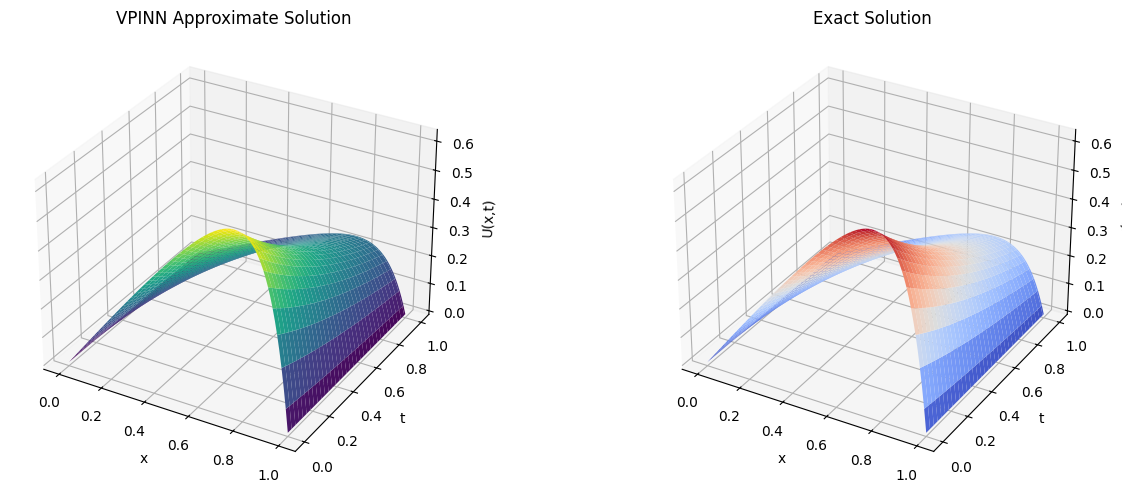}
	\caption{}
	\label{fig:subfig7A}
\end{subfigure}
\hfill
\begin{subfigure}[b]{1.0\textwidth}
	\includegraphics[width=\textwidth]{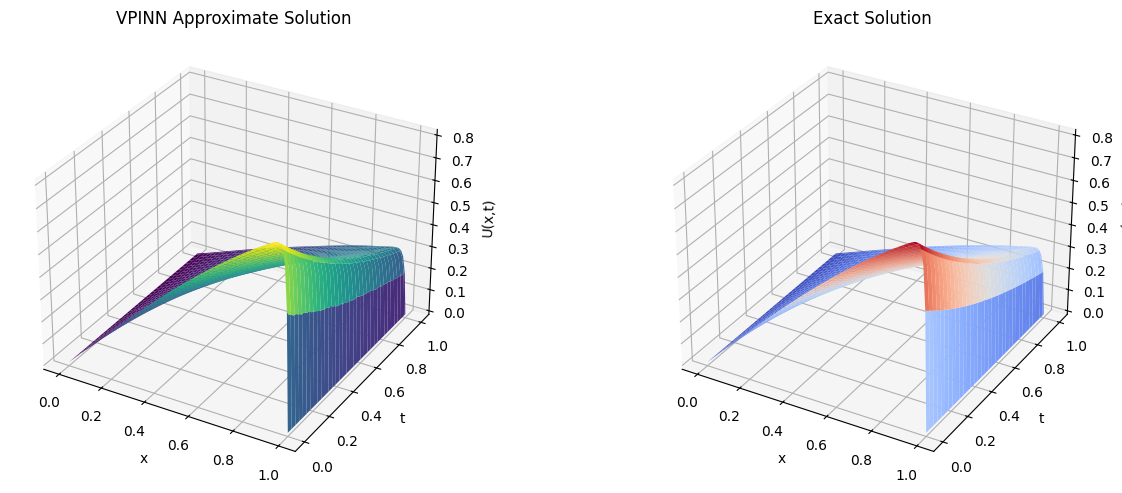}
	\caption{}
	\label{fig:subfig7B}
\end{subfigure}
\hfill
\begin{subfigure}[b]{1.0\textwidth}
	\includegraphics[width=\textwidth]{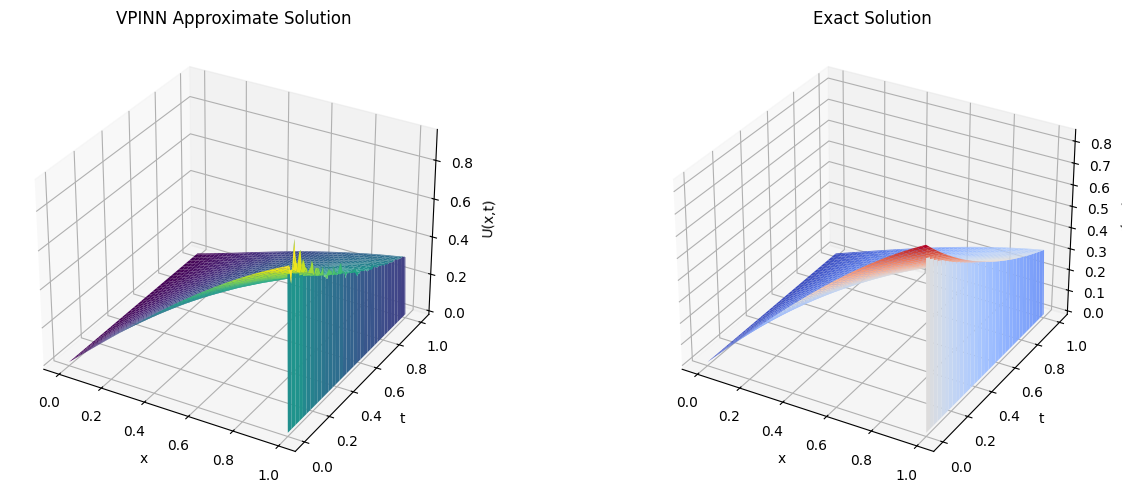}
	\caption{}
	\label{fig:subfig7C}
\end{subfigure}
\caption{Comparison between exact solution vs approximate solution for Example \ref{example2} with different value of perturbation parameter $\epsilon=10^{-1},10^{-2}$,and $10^{-3}$, corresponding to subfigures~(a), (b), and (c), respectively.}\label{figure7}
\end{figure}

We executed the code for $1000$ iterations, and the computed errors at the final time $T = 1$ are reported in Table~\ref{tbl_2}. Figure~\ref{figure4} illustrates the comparison between the exact solution and the VPINN approximation for various parameter values at $T = 1$. The corresponding absolute error distributions are shown in Figure~\ref{figure5}. Figure~\ref{figure6} displays the loss versus epoch curves, which represent the convergence behavior during the final iteration of the algorithm. Additionally, Figure~\ref{figure7} provides a surface plot comparison  between the exact solution and the VPINN approximation for various values of the perturbation parameter.

\begin{table}[htbp]
\caption{\it{  Loss and Error for Example \ref{example2} at the  final time $T=1$.}}\label{tbl_2}  
\begin{tabular}{@{}lllllll@{}}
	\multicolumn 1 {c}{}  & \multicolumn 5 {c}
	{}\\
	\hline
	$\epsilon$	& Loss & maximum & Rel-maximum  & $L_{2}$ & Rel-$L_{2}$ \\
	\hline
	$10^{-1}$  & 3.79781e-07 & 4.73222e-04 &  2.05343e-03 & 2.63204e-04 & 1.68000e-03\\ [6pt]
	
	$10^{-2}$  &  8.91816e-06 &4.82719e-03 & 1.62415e-02 & 1.33777e-03 & 7.13798e-03 \\ [6pt]
	
	$10^{-3}$  & 3.38240e-06 &  2.65777e-03 &  8.64250e-03 &  1.30504e-03 &  6.87290e-03 \\ [6pt]
	
	$10^{-4}$    &1.31400e-04 &  8.27751e-03  & 2.69155e-02&  2.06461e-03 &   1.08731e-02  \\ [6pt]
	
	$10^{-5}$    &  3.39702e-05 & 2.12925e-02  &6.92356e-02&   2.31294e-03 &1.21809e-02 \\ [6pt]
	
	\hline
\end{tabular}
\end{table}

\newpage
\subsection{The reaction-diffusion BVP}\label{sec3.3}
We examine the following one-parameter singularly perturbed BVPs:
\begin{equation}\label{E8}
\left\{
\begin{array} {ll}
	-\epsilon^{2} u''(x)+c(x)u(x)=r(x),\  x\in I=(0,\,1), \\[8pt]
	u(0)= u(1)=0.\\[8pt]
\end{array}\right.
\end{equation}
Here, $0<\epsilon\ll 1$ is a small  parameter, and the function  $c(x)$ and  $r(x)$ are sufficiently smooth and satisfy the conditions  $c(x)\geq \beta^{2} > 0$,  $\forall x \in I$, for constant $\beta$. Therefore solution $u(x)$ of (\ref{E8}) exhibits boundary layers  at $x=0$ and $x=1$.

The finite element formulation and algorithm will be similar to the model problem (\ref{E1}).


\begin{example}\label{example3}
Consider the problem in \eqref{E8} with \( c(x) = 1 \),
where the source function \( r(x) \) is computed from the exact solution
\[
u(x) = (1+exp(-1/\epsilon)-exp(-x/\epsilon)-exp(-(1-x)/\epsilon)).
\]
\end{example}

\begin{figure}[htbp]
\centering

\begin{subfigure}[b]{0.3\textwidth}
\includegraphics[width=\textwidth,height=1.2\textwidth]{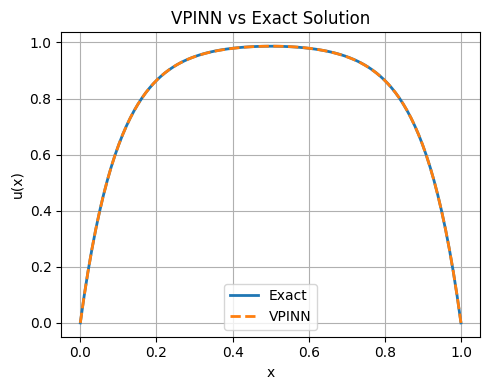}
\caption{}
\label{fig:subfig8A}
\end{subfigure}
\hfill
\begin{subfigure}[b]{0.3\textwidth}
\includegraphics[width=\textwidth,height=1.2\textwidth]{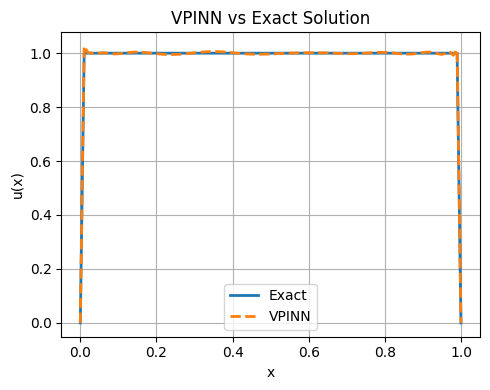}
\caption{}
\label{fig:subfig8B}
\end{subfigure}
\hfill
\begin{subfigure}[b]{0.3\textwidth}
\includegraphics[width=\textwidth,height=1.2\textwidth]{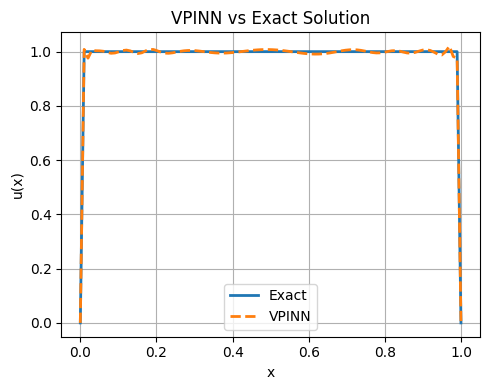}
\caption{}
\label{fig:subfig8C}
\end{subfigure}
\caption{Comparison between Exact and VPINN solution for Example \ref{example3} with different value of perturbation parameter $\epsilon=10^{-1},10^{-3}$, and $10^{-5}$, corresponding to subfigures~(a), (b), and (c), respectively.}\label{figure8}
\end{figure}

\begin{figure}[htbp]
\centering

\begin{subfigure}[b]{0.3\textwidth}
\includegraphics[width=\textwidth,height=1.0\textwidth]{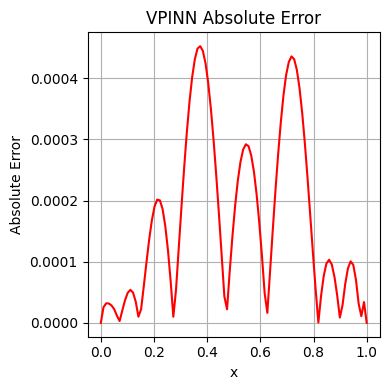}
\caption{}
\label{fig:subfig9A}
\end{subfigure}
\hfill
\begin{subfigure}[b]{0.3\textwidth}
\includegraphics[width=\textwidth,height=1.0\textwidth]{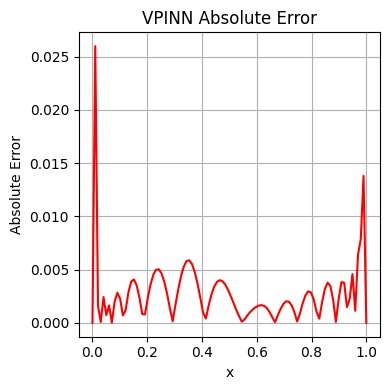}
\caption{}
\label{fig:subfig9B}
\end{subfigure}
\hfill
\begin{subfigure}[b]{0.3\textwidth}
\includegraphics[width=\textwidth,height=1.0\textwidth]{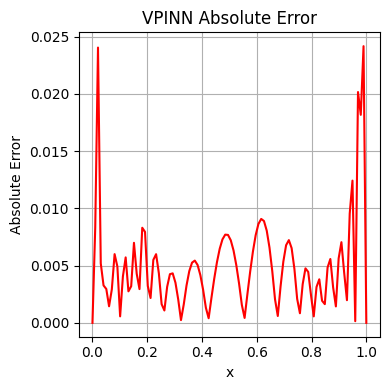}
\caption{}
\label{fig:subfig9C}
\end{subfigure}
\caption{Error vs x for Example \ref{example3} with different value of perturbation parameter$\epsilon=10^{-1},10^{-3},$ and $10^{-5}$, corresponding tosubfigures~(a), (b), and (c), respectively.}\label{figure9}
\end{figure}

\begin{figure}[htbp]
\centering

\begin{subfigure}[b]{0.3\textwidth}
\includegraphics[width=\textwidth,height=1.0\textwidth]{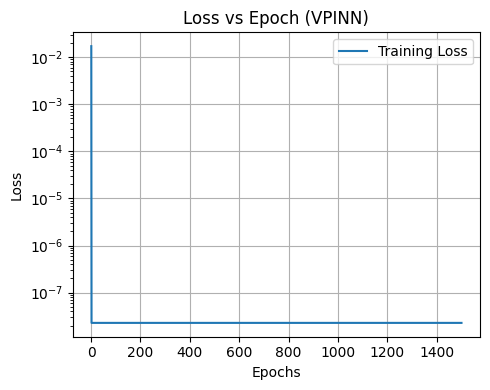}
\caption{}
\label{fig:subfig10A}
\end{subfigure}
\hfill
\begin{subfigure}[b]{0.3\textwidth}
\includegraphics[width=\textwidth,height=1.0\textwidth]{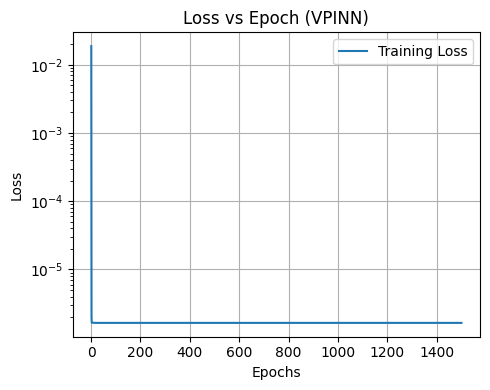}
\caption{}
\label{fig:subfig10B}
\end{subfigure}
\hfill
\begin{subfigure}[b]{0.3\textwidth}
\includegraphics[width=\textwidth,height=1.0\textwidth]{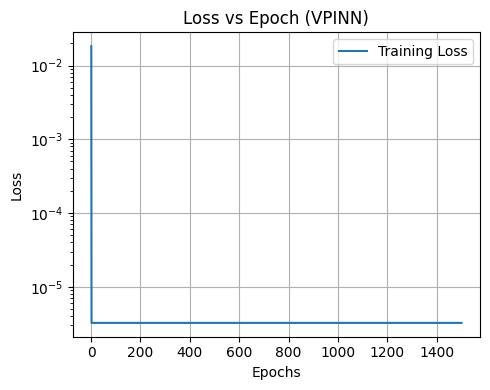}
\caption{}
\label{fig:subfig10C}
\end{subfigure}
\caption{Loss vs Epoch for Example \ref{example3} with different value of perturbation parameter $\epsilon=10^{-1},10^{-3},$ and $10^{-5}$, corresponding to subfigures~(a), (b), and (c), respectively.}\label{figure10}
\end{figure}

We executed the code for $1500$ iterations and the Maximum error, Relative maximum error, $L_{2}$ error, Relative $L_{2}$ error and Loss at $1500$ epoch obtained are presented in Table \ref{tbl_3}. Figure \ref{figure8} presents a comparison between the exact and VPINN solutions for different parameter values. Figure \ref{figure9} illustrates the absolute error plots for the same example, while Figure \ref{figure10} shows the loss versus epoch curves, corresponding to the final iteration of our algorithm.

\begin{table}[htbp]
\caption{\it{ Loss and Error for Example \ref{example3}.}}\label{tbl_3}  
\begin{tabular}{@{}lllllll@{}}
\multicolumn 1 {c}{}  & \multicolumn 5 {c}
{}\\
\hline
$\epsilon$	& Loss & maximum & Rel-maximum  & $L_{2}$ & Rel-$L_{2}$ \\
\hline
$10^{-1}$  &  2.97115e-08 & 4.52220e-04 & 4.58385e-04 & 2.16551e-04 & 2.60098e-05\\ [6pt]
$10^{-2}$  &  1.66846e-06 &  4.05634e-02  & 4.05634e-02 & 7.09868e-03 & 7.24451e-04\\ [6pt]

$10^{-3}$  & 2.29249e-06 &  2.59871e-02  &  2.59871e-02 &  4.08781e-03 &   4.12932e-04 \\ [6pt]

$10^{-4}$ & 5.23947e-06 &  3.31051e-02  & 3.31051e-02 &  6.11351e-03 &   6.17558e-04   \\ [6pt]

$10^{-5}$   & 2.53564e-06  & 2.41746e-02  & 2.41746e-02  &  6.53846e-03 & 6.60484e-04 \\ [6pt]

\hline
\end{tabular}
\end{table}

\subsection{Two-parameter SPPs}\label{sec3.4}
We examine the following two-parameter singularly perturbed BVPs:
\begin{equation}\label{E10}
\left\{
\begin{array} {ll}
-\epsilon u''(x)+ \mu b(x)u'(x)+c(x)u(x)=r(x),\  x\in I=(0,\,1), \\[8pt]
u(0)=u(1)=0.\\[8pt]
\end{array}\right.
\end{equation}
Here, $0<\epsilon\ll 1$ and   $0<\mu\ll 1$  are small perturbation parameter, and the function  $b(x), c(x)$ and  $r(x)$ are sufficiently smooth and satisfy the conditions  $b(x)\geq \alpha > 0$, \quad $c(x)-\frac{\mu b'(x)}{2}>0$,  $\forall x \in I$, for constant $\alpha$. Therefore solution $u(x)$ of (\ref{E10}) exhibits boundary layers both the ends of domain i.e, at $x=0$ and $x=1$. 

\begin{lemma}
Let $z_{0}$ and $z_{1}$ be the solution of the characteristic equation.
\[-\epsilon z^{2}(x)+\mu b(x)z(x)+c(x)=0.\]
Here,  $z_{0}(x) <0,\quad z_{1}(x)>0$ the boundary layers at $x=0$ and $x=1$, respectively.

$\lambda_{0}=-\max\limits_{x\in [0,1]} z_{0}(x)$, \quad  $\lambda_{1}=\min\limits_{x\in [0,1]} z_{1}(x)$,

\[\lambda_{0}=\max\limits_{x\in [0,1]} \frac{-\mu b(x)+\sqrt{\mu^{2}b^{2}(x)+4\epsilon c(x)}}{2\epsilon},\]

\[\lambda_{1}=\min\limits_{x\in [0,1]} \frac{\mu b(x)+\sqrt{\mu^{2}b^{2}(x)+4\epsilon c(x)}}{2\epsilon}.\]
\end{lemma}

Depending on the relationship between $\epsilon$ and $\mu$ , we can classify (\ref{E10}) into the three regimes:
\begin{itemize}
\item If $\epsilon \ll \mu = 1 $, then $\lambda_{0} = O(1)$ and $\lambda_{1} = O(1/\epsilon)$. In this case is similar to convection-diffusion type.

\item If $\epsilon \ll \mu^{2} \ll 1 $, then $\lambda_{0} = O(\mu^{-1})$ and $\lambda_{1} = O(\mu\epsilon^{-1})$. In this case is similar to diffusion-convection-reaction type.

\item If $ \mu^{2} \ll \epsilon \ll 1 $, then $\lambda_{0} = O(\mu^{-1/2})$ and $\lambda_{1} = O(\epsilon^{-1/2})$. In this case is similar to reaction-diffusion type.
\end{itemize}

%

\begin{example}\label{example4}
Consider the problem in \eqref{E10} with 
\( b(x) = 1 \) and \( c(x) = 1 \).
The source function \( r(x) \) is chosen from the exact solution
\[
u(x)
= \frac{1}{E}\left[
\frac{ e^{m_1}\big(e^{1-m_1}-1\big)\,e^{-m_2 x}
- \big(e^{1-m_2}-1\big)\,e^{m_1(1-x)} }{D}
- e^{\,1-x}
\right],
\]
where 
\(
m_{1,2} = \frac{\mu \mp \sqrt{\mu^{2} + 4 \epsilon}}{2 \epsilon},\quad
D = 1 - \exp\!\left(\frac{-\sqrt{\mu^{2} + 4 \epsilon}}{\epsilon}\right),\quad
E = \epsilon - \mu - 1.
\)
\end{example}

\begin{figure}[htbp]
\centering

\begin{subfigure}[b]{0.3\textwidth}
\includegraphics[width=\textwidth,height=1.2\textwidth]{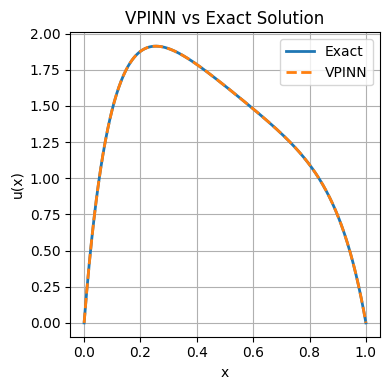}
\caption{$(\epsilon,\mu)=(10^{-2},10^{-3})$}
\label{fig:subfig11A}
\end{subfigure}
\hfill
\begin{subfigure}[b]{0.3\textwidth}
\includegraphics[width=\textwidth,height=1.2\textwidth]{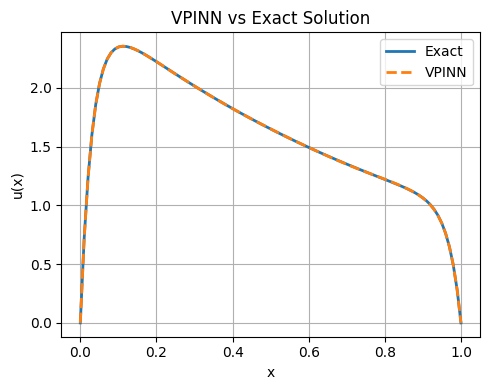}
\caption{$(\epsilon,\mu)=(10^{-3},10^{-4})$}
\label{fig:subfig11B}
\end{subfigure}
\hfill
\begin{subfigure}[b]{0.3\textwidth}
\includegraphics[width=\textwidth,height=1.2\textwidth]{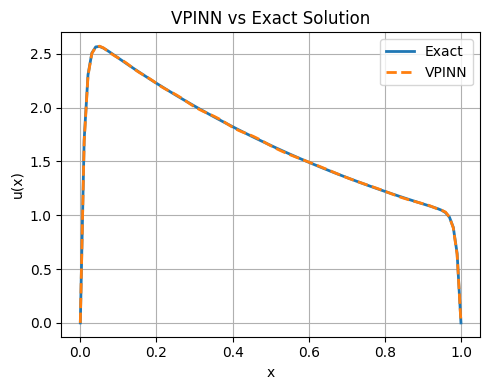}
\caption{$(\epsilon,\mu)=(10^{-4},10^{-5})$}
\label{fig:subfig11C}
\end{subfigure}
\caption{Comparison between Exact and VPINN solution for Example \ref{example4} with different value of perturbation parameter.}\label{figure11}
\end{figure}

\begin{figure}[htbp]
\centering

\begin{subfigure}[b]{0.3\textwidth}
\includegraphics[width=\textwidth,height=1.0\textwidth]{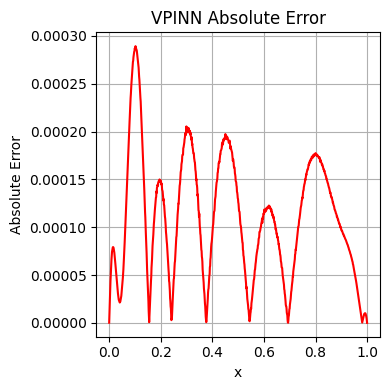}
\caption{$(\epsilon,\mu)=(10^{-2},10^{-3})$}
\label{fig:subfig12A}
\end{subfigure}
\hfill
\begin{subfigure}[b]{0.3\textwidth}
\includegraphics[width=\textwidth,height=1.0\textwidth]{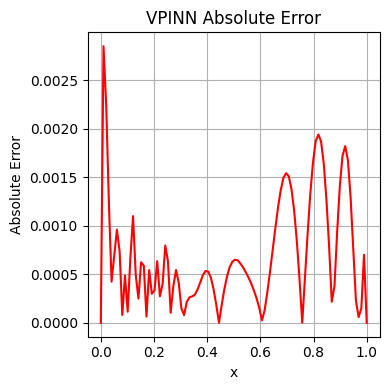}
\caption{$(\epsilon,\mu)=(10^{-3},10^{-4})$}
\label{fig:subfig12B}
\end{subfigure}
\hfill
\begin{subfigure}[b]{0.3\textwidth}
\includegraphics[width=\textwidth,height=1.0\textwidth]{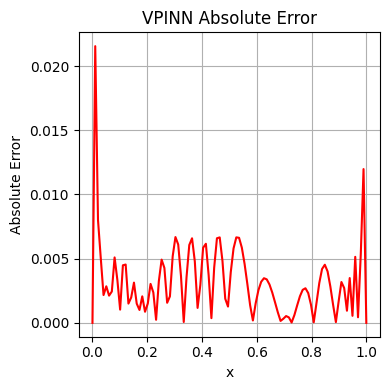}
\caption{$(\epsilon,\mu)=(10^{-4},10^{-5})$}
\label{fig:subfig12C}
\end{subfigure}
\caption{Error vs x for Example \ref{example4} with different value of perturbation parameter.}\label{figure12}
\end{figure}

\begin{figure}[htbp]
\centering

\begin{subfigure}[b]{0.3\textwidth}
\includegraphics[width=\textwidth,height=1.0\textwidth]{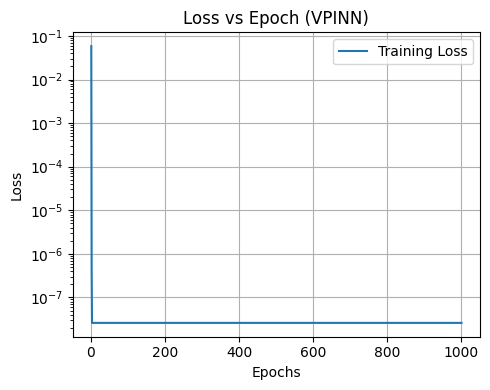}
\caption{$(\epsilon,\mu)=(10^{-2},10^{-3})$}
\label{fig:subfig13A}
\end{subfigure}
\hfill
\begin{subfigure}[b]{0.3\textwidth}
\includegraphics[width=\textwidth,height=1.0\textwidth]{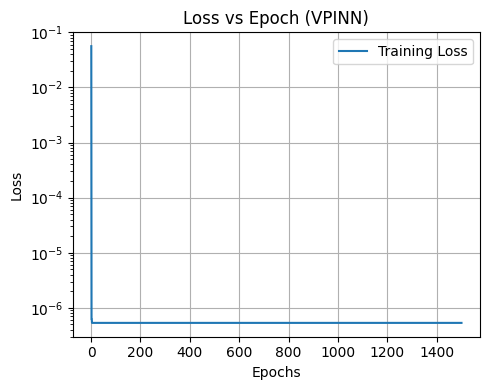}
\caption{$(\epsilon,\mu)=(10^{-3},10^{-4})$}
\label{fig:subfig13B}
\end{subfigure}
\hfill
\begin{subfigure}[b]{0.3\textwidth}
\includegraphics[width=\textwidth,height=1.0\textwidth]{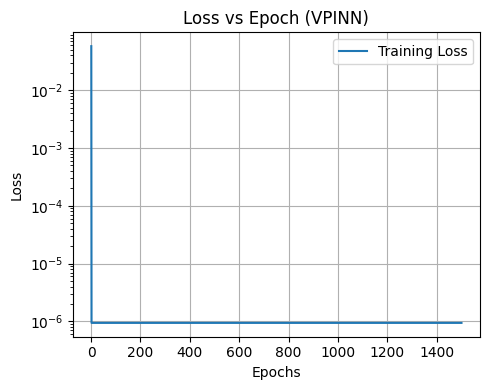}
\caption{$(\epsilon,\mu)=(10^{-4},10^{-5})$}
\label{fig:subfig13C}
\end{subfigure}
\caption{Loss vs Epoch for Example \ref{example4} with different value of perturbation parameter.}\label{figure13}
\end{figure}

We executed the code for $1500$ iterations and the Maximum error, Relative maximum error, $L_{2}$ error, Relative $L_{2}$ error and Loss at $1500$ epoch obtained are presented in Table \ref{tbl_4}. Figure \ref{figure11} presents a comparison between the exact and VPINN solutions for different parameter values. Figure \ref{figure12} illustrates the absolute error plots for the same example, while Figure \ref{figure13} shows the loss versus epoch curves, corresponding to the final iteration of our algorithm.

\begin{table}[htbp]
\caption{\it{  Loss and Error for Example \ref{example4}.}}\label{tbl_4}  
\begin{tabular}{@{}lllllll@{}}
\multicolumn 1 {c}{}  & \multicolumn 5 {c}
{}\\
\hline
$\epsilon,\ \ \quad \mu$	& Loss & maximum & Rel-maximum  & $L_{2}$ & Rel-$L_{2}$ \\
\hline
$10^{-1},10^{-2}$  &  3.06321e-07 & 3.80754e-04  & 3.66852e-04 & 2.03801e-04 & 8.37204e-06 \\ [6pt]
$10^{-2},10^{-3}$  &  2.61157e-08 & 2.89321e-04 & 1.51199e-04 & 1.27626e-04 & 2.79353e-06 \\ [6pt]

$10^{-3},10^{-4}$ &  5.32722e-07& 2.85333e-03  & 1.21211e-03 &  8.97248e-04 &  5.37536e-04  \\ [6pt]

$10^{-4},10^{-5}$ & 9.44896e-07  & 2.15509e-02  &  8.39448e-03 &  4.27543e-03 &   2.45244e-03 \\ [6pt]
\hline
\end{tabular}
\end{table}

\subsection{Parabolic two-parameter IBVP }\label{sec3.5}
Consider the following singularly perturbed parabolic PDEs with
two-parameters in the domain $G$ = $I$ $\times$ $(0, T]$, where $I = (0, 1):$
\begin{equation}\label{E12}
\left\{
\begin{array} {ll}
u_{t}	-\epsilon u_{xx}+\mu b(x)u_{x}+c(x)u(x,t)=r(x,t),\  x\in I=(0,\,1), \\[8pt]
u(0,t)= u(1,t)=0,  \quad 0\leq t \leq T,\\[8pt]
u(x,0)-\phi(x)=0, \quad \mbox{for}\quad  x \in [0, 1].
\end{array}\right.
\end{equation}
Here, $0<\epsilon\ll 1$ and $0<\mu\ll 1$  are small perturbation parameter, and the function  $b(x), c(x)$ and  $r(x,t)$ are sufficiently smooth and satisfy the conditions  $b(x)\geq \alpha > 0$, \quad $c(x)-\frac{\mu b'(x)}{2}>0$,  $\forall x \in I$, for constant $\alpha$. Therefore solution $u(x,t)$ of (\ref{E12}) exhibits boundary layers of both the ends of domain i.e, at $x=0$ and $x=1$.

\begin{figure}[htbp]
\centering

\begin{subfigure}[b]{0.3\textwidth}
\includegraphics[width=\textwidth,height=1.2\textwidth]{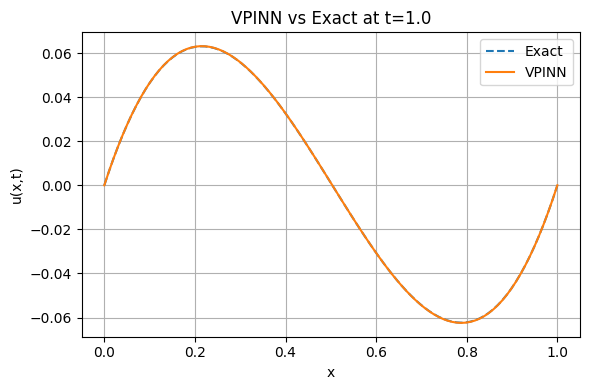}
\caption{$(\epsilon,\mu)=(10^{-1},10^{-2})$}
\label{fig:subfig14A}
\end{subfigure}
\hfill
\begin{subfigure}[b]{0.3\textwidth}
\includegraphics[width=\textwidth,height=1.2\textwidth]{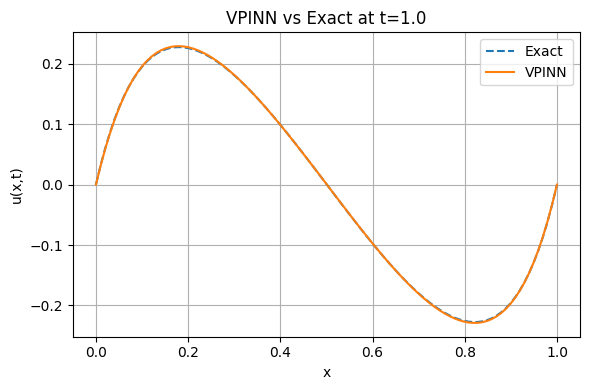}
\caption{$(\epsilon,\mu)=(10^{-2},10^{-3})$}
\label{fig:subfig14B}
\end{subfigure}
\hfill
\begin{subfigure}[b]{0.3\textwidth}
\includegraphics[width=\textwidth,height=1.2\textwidth]{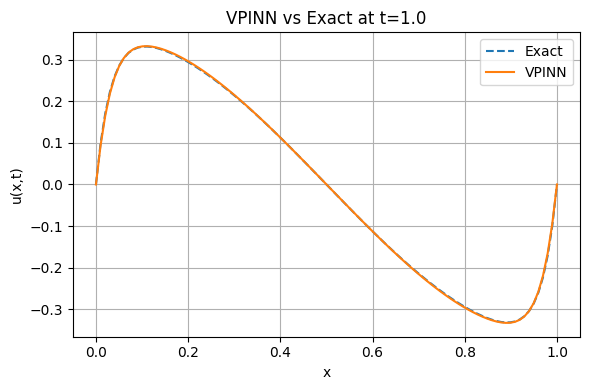}
\caption{$(\epsilon,\mu)=(10^{-3},10^{-4})$}
\label{fig:subfig14C}
\end{subfigure}
\caption{Comparison between Exact and VPINN solution for Example \ref{example5} with different value of perturbation parameter.}\label{figure14}
\end{figure}
\begin{figure}[htbp]
\centering

\begin{subfigure}[b]{0.3\textwidth}
\includegraphics[width=\textwidth,height=1.0\textwidth]{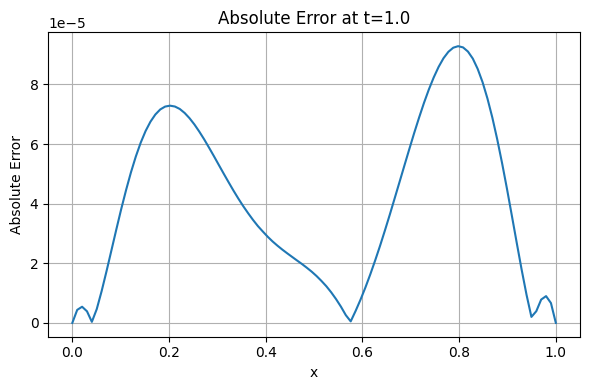}
\caption{$(\epsilon,\mu)=(10^{-1},10^{-2})$}
\label{fig:subfig15A}
\end{subfigure}
\hfill
\begin{subfigure}[b]{0.3\textwidth}
\includegraphics[width=\textwidth,height=1.0\textwidth]{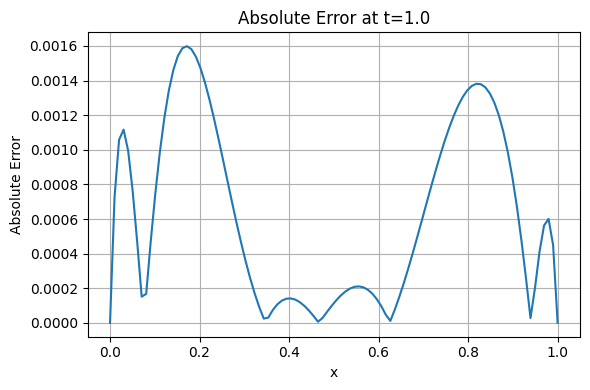}
\caption{$(\epsilon,\mu)=(10^{-2},10^{-3})$}
\label{fig:subfig15B}
\end{subfigure}
\hfill
\begin{subfigure}[b]{0.3\textwidth}
\includegraphics[width=\textwidth,height=1.0\textwidth]{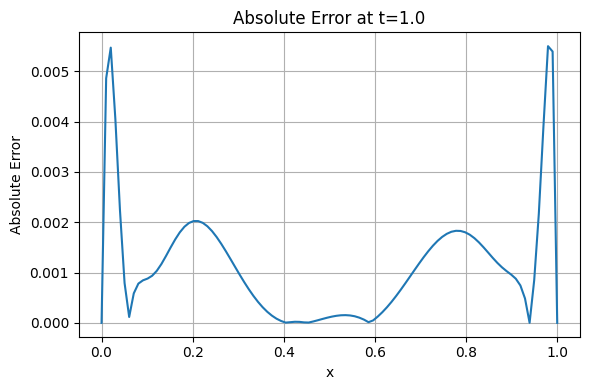}
\caption{$(\epsilon,\mu)=(10^{-3},10^{-4})$}
\label{fig:subfig15C}
\end{subfigure}
\caption{Error vs x  for Example \ref{example5} with different value of perturbation parameter.}\label{figure15}
\end{figure}
\begin{figure}[htbp]
\centering

\begin{subfigure}[b]{0.3\textwidth}
\includegraphics[width=\textwidth,height=1.0\textwidth]{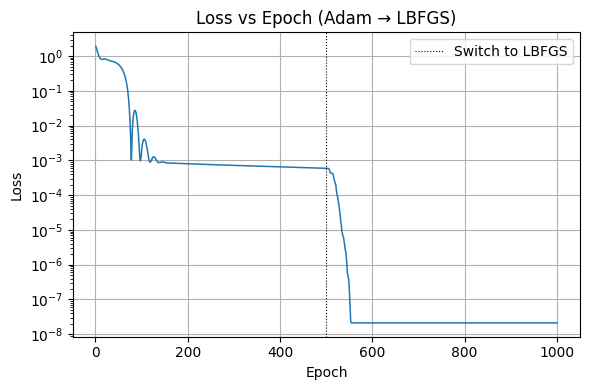}
\caption{$(\epsilon,\mu)=(10^{-1},10^{-2})$}
\label{fig:subfig16A}
\end{subfigure}
\hfill
\begin{subfigure}[b]{0.3\textwidth}
\includegraphics[width=\textwidth,height=1.0\textwidth]{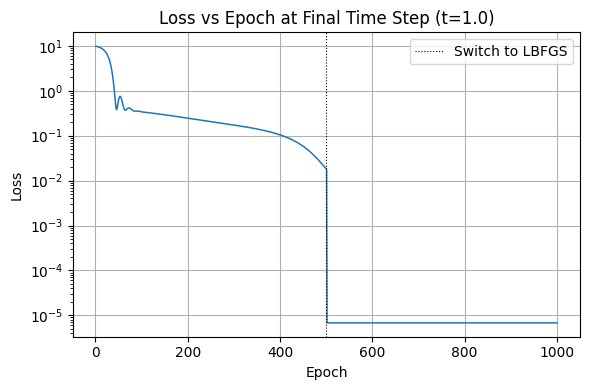}
\caption{$(\epsilon,\mu)=(10^{-2},10^{-3})$}
\label{fig:subfig16B}
\end{subfigure}
\hfill
\begin{subfigure}[b]{0.3\textwidth}
\includegraphics[width=\textwidth,height=1.0\textwidth]{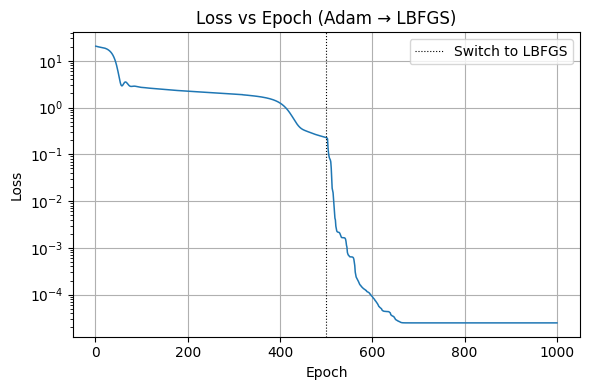}
\caption{$(\epsilon,\mu)=(10^{-3},10^{-4})$}
\label{fig:subfig16C}
\end{subfigure}
\caption{Loss vs Epoch for Example \ref{example5} with different value of perturbation parameter.}\label{figure16}
\end{figure}
\begin{figure}[htbp]
\centering
\begin{subfigure}[b]{1.0\textwidth}
\includegraphics[width=\textwidth]{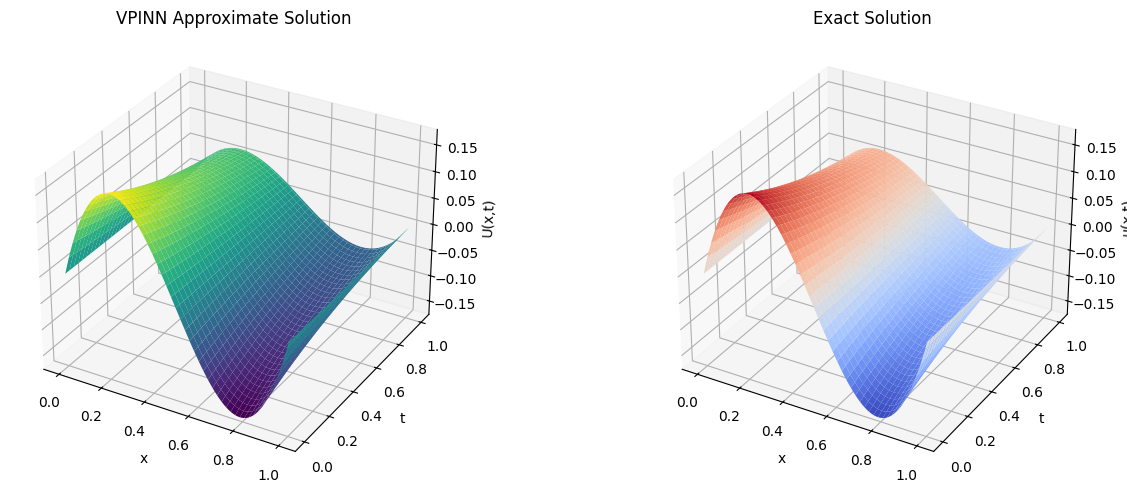}
\caption{$(\epsilon,\mu)=(10^{-1},10^{-2})$}
\label{fig:subfig17A}
\end{subfigure}
\hfill
\begin{subfigure}[b]{1.0\textwidth}
\includegraphics[width=\textwidth]{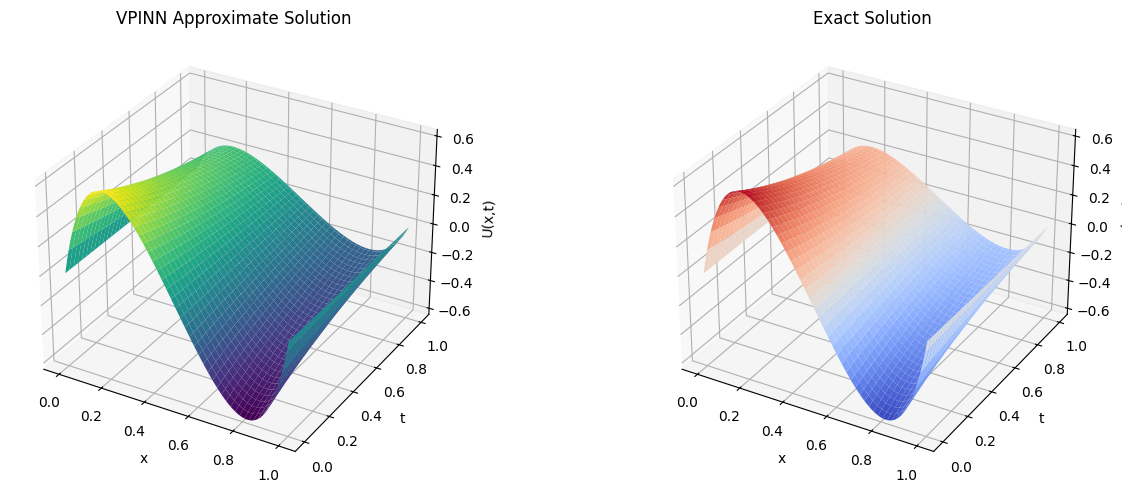}
\caption{$(\epsilon,\mu)=(10^{-2},10^{-3})$}
\label{fig:subfig17B}
\end{subfigure}
\hfill
\begin{subfigure}[b]{1.0\textwidth}
\includegraphics[width=\textwidth]{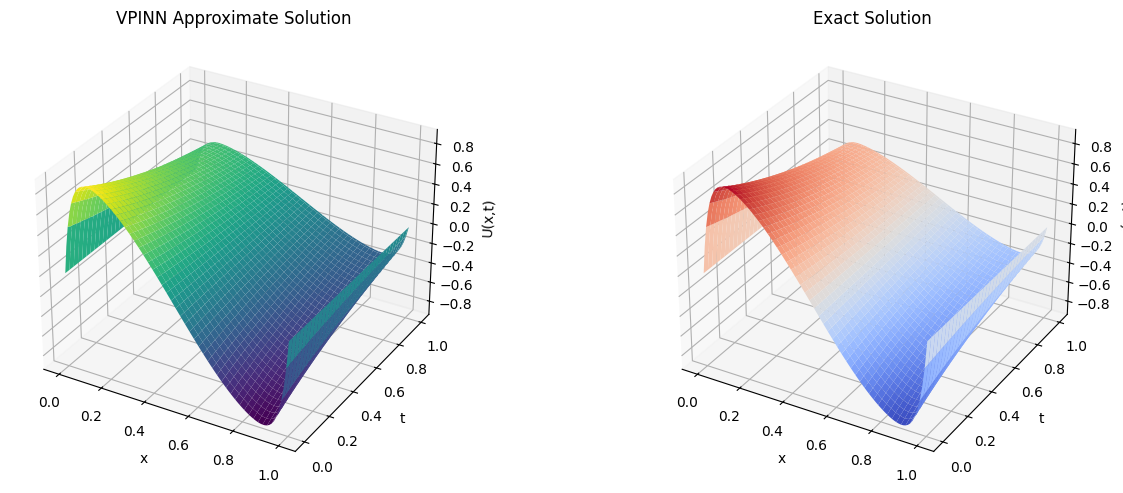}
\caption{$(\epsilon,\mu)=(10^{-3},10^{-4})$}
\label{fig:subfig17C}
\end{subfigure}
\caption{Comparison between exact solution vs approximate solution for Example \ref{example5} with perturbation parameter.}\label{figure17}
\end{figure}

\begin{example}\label{example5}
Consider  the parabolic two-parameter IBVP \eqref{E12} with
\( b(x) = 1 \) and \( c(x) = 1 \).
The source term \( r(x,t) \) and initial condition \( \phi(x) \) 
are chosen from the exact solution
\[
u(x,t) = e^{-t}\bigg(a \cos(\pi x) + b \sin(\pi  x) + A \exp(-u_{l} x) +B \exp\big(-u_{r} (1 - x)\big)\bigg),\]

where \quad $a = \frac{(\epsilon \pi^{2} + 1)}{\mu^{2} \pi^{2} + (\epsilon \pi^{2} + 1)^{2}}$,\quad $b = \frac{\mu \pi}{\mu^{2} \pi^{2} + (\epsilon \pi^{2} + 1)^{2}}$,\quad
$u_{l,r} = \frac{(\mp\mu + \sqrt{\mu^{2} + 4 \epsilon})}{2 \epsilon}$,\quad
$A =-a\frac{ 1 + exp(-u_r)}{1 - exp\big(-(u_l+u_r)\big)}$,\quad
$B = a\frac{ 1 + exp(-u_l)}{1 - exp\big(-(u_l+u_r)\big)}$.
\end{example}
We executed the code for $1000$ iterations, and the computed errors at the final time $T = 1$ are reported in Table~\ref{tbl_5}. Figure~\ref{figure14} illustrates the comparison between the exact solution and the VPINN approximation for various parameter values at $T = 1$. The corresponding absolute error distributions are shown in Figure~\ref{figure15}. Figure~\ref{figure16} displays the loss versus epoch curves, which represent the convergence behavior during the final iteration of the algorithm. Additionally, Figure~\ref{figure17} provides a surface plot comparison between the exact solution and the VPINN approximation for various values of the perturbation parameter.

\begin{table}[htbp]
\caption{\it{  Loss and Error for Example \ref{example5} at the final time $T=1$.}}\label{tbl_5}  
\begin{tabular}{@{}lllllll@{}}
\multicolumn 1 {c}{}  & \multicolumn 5 {c}
{}\\
\hline
$\epsilon,\ \ \quad \mu$	 & Loss & maximum & Rel-maximum  & $L_{2}$ & Rel-$L_{2}$\\
\hline
$10^{-1},10^{-2}$  &  2.11662e-08  & 9.28079e-05  &  1.47207e-03 &   4.98550e-05 & 1.11424e-03 \\ [6pt]

$10^{-2},10^{-3}$ &  6.82480e-06  &  1.59839e-03  &  7.02785e-03 & 8.08532e-04 & 4.99766e-03 \\ [6pt]

$10^{-3},10^{-4}$ & 2.49291e-05  &  5.49714e-03 &  1.65834e-02 &  1.62604e-03 &  7.03748e-03\\ [6pt]

$10^{-4},10^{-5}$ & 3.25968e-05 &  2.22948e-02 &  6.18356e-02 &  3.18900e-03 &   1.27222e-02   \\ [6pt]

\hline
\end{tabular}
\end{table}

\section{Conclusion}\label{Conclusion}

In this article, we presented a VPINN framework based on the Petrov-Galerkin method for solving SPPs. The solution was approximated using DNNs for the trial space, while the test space consisted of localized hat functions. This setup allowed for effective handling of sharp boundary layers without relying on dense mesh discretizations. For time-dependent problems, we employed the backward Euler scheme for temporal discretization. The training process combined the strengths of both ADAM and L-BFGS optimizers starting with ADAM for faster initial progress and then refining with L-BFGS. Source terms were evaluated through automatic differentiation, enhancing accuracy within the weak formulation. The proposed VPINN method demonstrated strong performance, offering a stable, accurate, and mesh-free solution technique for both steady and time-dependent SPPs.

\section*{Funding}
The authors received no external funding for this work.

\section*{Conflict of Interest}
The authors declare that they have no conflicts of interest.

\section*{Data Availability}
Not applicable.


\clearpage
\begin{thebibliography}{00}

\bibitem{Alhu}
K.~Alhumaizi.
\newblock Flux-limiting solution techniques for simulation of
reaction-diffusion-convection system.
\newblock {\em Commun. Nonlinear Sci. Numer. Simul.}, 12(6):953--965, 2007.

\bibitem{IPDG}
D.~N.~Arnold.
\newblock An interior penalty finite element method with discontinuous
elements.
\newblock {\em SIAM journal on numerical analysis}, 19(4):742--760, 1982.

\bibitem{Parabolic_FEM_2D}
D.~Avijit and S.~Natesan.
\newblock A novel two-step streamline-diffusion {FEM} for singularly perturbed
2{D} parabolic {PDE}s.
\newblock {\em Appl. Numer. Math.}, 172:259--278, 2022.

\bibitem{PINN_Natesan}
P.~Boro, A.~Raina, and S.~Natesan.
\newblock A parameter-driven physics-informed neural network framework for
solving two-parameter singular perturbation problems involving boundary
layers, 2025. https://arxiv.org/abs/2505.01159


\bibitem{PAPINN}
F.~Cao, F,~Gao, X.~Guo, and D.~Yuan.
\newblock Physics-informed neural networks with parameter asymptotic strategy
for learning singularly perturbed convection-dominated problem.
\newblock {\em Comput. Math. Appl.}, 150:229--242, 2023.

\bibitem{DNN_FEM_4}
G.~Capuano and J.~J.~Rimoli.
\newblock Smart finite elements: a novel machine learning application.
\newblock {\em Comput. Methods Appl. Mech. Engrg.}, 345:363--381, 2019.

\bibitem{PINN_NIPG}
G.~Chen, S.~Xu, D.~Ni, and T.~Zeng.
\newblock Dgnn: A neural pde solver induced by discontinuous galerkin methods,
2025. https://arxiv.org/abs/2503.10021

\bibitem{SIPG}
R.~Hartmann and P.~Houston.
\newblock Symmetric interior penalty {DG} methods for the compressible
{N}avier-{S}tokes equations. {II}. {G}oal-oriented a posteriori error
estimation.
\newblock {\em Int. J. Numer. Anal. Model.}, 3(2):141--162, 2006.

\bibitem{DNN_FEM_3}
J.~S.~Hesthaven and S.~Ubbiali.
\newblock Non-intrusive reduced order modeling of nonlinear problems using
neural networks.
\newblock {\em J. Comput. Phys.}, 363:55--78, 2018.

\bibitem{VPINN1}
E.~Kharazmi, Z.~Zhang, and G.~E. Karniadakis.
\newblock Variational physics-informed neural networks for solving partial
differential equations, 2019. https://arxiv.org/abs/1912.00873

\bibitem{VPINN}
E.~Kharazmi, Z.~Zhang, and G.~E. Karniadakis.
\newblock {$hp$}-{VPINN}s: variational physics-informed neural networks with
domain decomposition.
\newblock {\em Comput. Methods Appl. Mech. Engrg.}, 374:113547, 25,
2021.




\bibitem{kokotovic}
P.~V.~Kokotovi\'c.
\newblock Applications of singular perturbation techniques to control problems.
\newblock {\em SIAM Rev.}, 26(4):501--550, 1984.

\bibitem{Modeling_Fluid}
Heinz-Otto Kreiss and Jens Lorenz.
\newblock {\em Initial-boundary value problems and the {N}avier-{S}tokes
equations}, volume~47 of {\em Classics in Applied Mathematics}.
\newblock Society for Industrial and Applied Mathematics (SIAM), Philadelphia,
PA, 2004.
\newblock Reprint of the 1989 edition.

\bibitem{DNN_FEM_1}
I.E. Lagaris, A.~Likas, and D.I. Fotiadis.
\newblock Artificial neural networks for solving ordinary and partial
differential equations.
\newblock {\em IEEE Transactions on Neural Networks}, 9(5):987--1000, 1998.

\bibitem{Deep_decomposition}
Ke~Li, Kejun Tang, Tianfan Wu, and Qifeng Liao.
\newblock D3m: A deep domain decomposition method for partial differential
equations.
\newblock {\em IEEE Access}, 8:5283--5294, 2020.



\bibitem{Melenk}
J.~M. Melenk and C.~Xenophontos.
\newblock Robust exponential convergence of {$hp$}-{FEM} in balanced norms for
singularly perturbed reaction-diffusion equations.
\newblock {\em Calcolo}, 53(1):105--132, 2016.

\bibitem{FB_PINN}
B.~Moseley, A.~Markham, and T.~N.~Meyer.
\newblock Finite basis physics-informed neural networks ({FBPINN}s): a scalable
domain decomposition approach for solving differential equations.
\newblock {\em Adv. Comput. Math.}, 49(4): 62, 39, 2023.



\bibitem{FB-PINN_Natesan}
A.~Raina, S.~Badireddi, and S.~Natesan.
\newblock Application of pinn to obtain solution of boundary layer problems
arising in fluid dynamics.
\newblock {\em Mathematical Foundations of Computing}, 2025. 10.3934/mfc.2025024

\bibitem{PINN_Raissi}
M.~Raissi, P.~Perdikaris, and G.~E. Karniadakis.
\newblock Physics-informed neural networks: a deep learning framework for
solving forward and inverse problems involving nonlinear partial differential
equations.
\newblock {\em J. Comput. Phys.}, 378:686--707, 2019.



\bibitem{Water_pollution}
A.~Rap, L.~Elliott, D.~B. Ingham, D.~Lesnic, and X.~Wen.
\newblock The inverse source problem for the variable coefficients
convection-diffusion equation.
\newblock {\em Inverse Probl. Sci. Eng.}, 15(5):413--440, 2007.

\bibitem{Reed}
W.~H.~Reed and T.~R.~Hill.
\newblock Triangular mesh methods for the neutron transport equation.
\newblock Technical report, Los Alamos Scientific Lab., N. Mex.(USA), 1973.

\bibitem{Riviere}
B.~Rivi\`ere.
\newblock {\em Discontinuous {G}alerkin methods for solving elliptic and
parabolic equations}, volume~35 of {\em Frontiers in Applied Mathematics}.
\newblock Society for Industrial and Applied Mathematics (SIAM), Philadelphia,
PA, 2008.
\newblock Theory and implementation.

\bibitem{Roos2}
H.~G.~Roos and M.~Schopf.
\newblock Convergence and stability in balanced norms of finite element methods
on {S}hishkin meshes for reaction-diffusion problems.
\newblock {\em ZAMM Z. Angew. Math. Mech.}, 95(6):551--565, 2015.

\bibitem{Schlicting}
H.~Schlichting and K.~Gersten.
\newblock {\em Boundary-layer theory}.
\newblock Springer-Verlag, Berlin, ninth edition, 2017.

\bibitem{Gautam1}
G.~Singh and S.~Natesan.
\newblock Study of the {NIPG} method for two-parameter singular perturbation
problems on several layer adapted grids.
\newblock {\em J. Appl. Math. Comput.}, 63(1-2):683--705, 2020.

\bibitem{Gautam_parabolic}
G.~Singh and S.~Natesan.
\newblock Superconvergence error estimates of discontinuous {G}alerkin time
stepping for singularly perturbed parabolic problems.
\newblock {\em Numer. Algorithms}, 90(3):1073--1090, 2022.

\bibitem{Gautam4}
M.~K.~Singh, G.~Singh, and S.~Natesan.
\newblock A unified study on superconvergence analysis of {G}alerkin {FEM} for
singularly perturbed systems of multiscale nature.
\newblock {\em J. Appl. Math. Comput.}, 66(1-2):221--243, 2021.

\bibitem{DNN_FEM_2}
J.~Sirignano and K.~Spiliopoulos.
\newblock D{GM}: a deep learning algorithm for solving partial differential
equations.
\newblock {\em J. Comput. Phys.}, 375:1339--1364, 2018.

\bibitem{fem_tobiska}
L.~Tobiska.
\newblock Analysis of a new stabilized higher order finite element method for
advection-diffusion equations.
\newblock {\em Comput. Methods Appl. Mech. Engrg.}, 196(1-3):538--550, 2006.

\bibitem{Yadav2}
S.~Yadav and S.~Ganesan.
\newblock Artificial neural network-augmented stabilized finite element method.
\newblock {\em J. Comput. Phys.}, 499: 112702, 15, 2024.

\bibitem{FEM_2D}
Z.~Zhang.
\newblock Finite element superconvergence on {S}hishkin mesh for 2-{D}
convection-diffusion problems.
\newblock {\em Math. Comp.}, 72(243):1147--1177, 2003.



\end{thebibliography}
\end{document}